\documentclass[onefignum,onetabnum,reqno]{siamart190516}

\usepackage{makecell}
\usepackage[utf8]{inputenc} 
\usepackage[T1]{fontenc}    
\usepackage{hyperref}       
\usepackage{url}            
\usepackage{booktabs}       
\usepackage{amsfonts}       
\usepackage{nicefrac}       
\usepackage{microtype}      
\usepackage{bbold}
\usepackage{capt-of}
\usepackage{mathtools}
\usepackage{mathabx}
\usepackage[center]{caption}
\usepackage{tablefootnote}
\usepackage{float}

\usepackage{amsmath, amssymb, graphicx, cite, epsfig,epstopdf,color,soul}
\usepackage{algorithmic}
\usepackage{tabularx}
\allowdisplaybreaks
\usepackage{color}
\usepackage{enumerate}
\usepackage[shortlabels]{enumitem}
\usepackage{multicol}
\usepackage{empheq}
\usepackage{caption} 
\usepackage{ulem}

\definecolor{blue}{gray}{0}

\newcommand{\Eset}{\mathbb{E}}

\newcommand{\Rset}{\mathbb{R}}

\newcommand{\Acal}{{\cal A}}

\newcommand{\Dcal}{{\cal D}}

\newcommand{\Fcal}{{\cal F}}

\newcommand{\Kcal}{{\cal K}}

\newcommand{\Mcal}{{\cal M}}

\newcommand{\Ocal}{{\cal O}}
\newcommand{\Pcal}{{\cal P}}

\newcommand{\Scal}{{\cal S}}

\newcommand{\Xcal}{{\cal X}}

\newcommand{\1}{{\mathbf{1}}}

\newcommand{\argmin}{\mathop{\rm argmin}}
\newcommand{\argmax}{\mathop{\rm argmax}}

\newtheorem{prop}{Proposition}
\newtheorem{lem}{Lemma}
\newtheorem{thm}{Theorem}
\newtheorem{cor}{Corollary}

\newtheorem{assump}{Assumption}

\usepackage{tikz}

\title{A Two-Time-Scale Stochastic Optimization Framework with Applications in Control and Reinforcement Learning}
\author{Sihan Zeng\thanks{School of Electrical and Computer Engineering, Georgia Institute of Technology, Atlanta, GA
.}
\and Thinh T. Doan\thanks{Bradley Department of Electrical and Computer Engineering, Virginia Tech, Blacksburg, VA.}
\and Justin Romberg\footnotemark[1]}

\begin{document}

\maketitle

\begin{abstract}
We study a new two-time-scale stochastic gradient method for solving optimization problems,
where the gradients are computed with the aid of an auxiliary variable under samples generated by time-varying Markov random processes controlled by the underlying optimization variable. These time-varying samples make gradient directions in our update biased and dependent, which can potentially lead to the divergence of the iterates. In our two-time-scale approach, one scale is to estimate the true gradient from these samples, which is then used to update the estimate of the optimal solution.
While these two iterates are implemented simultaneously, the former is updated ``faster'' (using bigger step sizes) than the latter (using smaller step sizes). Our first contribution is to characterize the finite-time complexity of the proposed two-time-scale stochastic gradient method. In particular, we provide explicit formulas for the convergence rates of this method under different structural assumptions, namely, strong convexity, Polyak-\L ojasiewicz (P\L) condition, and general non-convexity.

We apply our framework to policy optimization problems in control and reinforcement learning. 
First, we look at the infinite-horizon average-reward Markov decision process with finite state and action spaces and derive a convergence rate of $\widetilde{\Ocal}(k^{-2/5})$ for the online actor-critic algorithm under function approximation, which recovers the best known rate derived specifically for this problem. Second, we study the linear-quadratic regulator and show that an online actor-critic method converges with rate $\widetilde{\Ocal}(k^{-2/3})$. Third, we use the actor-critic algorithm to solve the policy optimization problem in an entropy regularized Markov decision process, where we also establish a convergence of $\widetilde{\Ocal}(k^{-2/3})$.
The results we derive for both the second and third problem are novel and previously unknown in the literature. Finally, we briefly present the application of our framework to gradient-based policy evaluation algorithms in reinforcement learning.
\end{abstract}

  


\section{Introduction}\label{sec:introduction}

Actor-critic algorithms are an important class of data-driven techniques for policy optimization in reinforcement learning.  They can be re-cast as optimization programs with a specific type of stochastic oracle for gradient evaluations.  In this paper, we present an abstraction of the actor-critic framework for solving general optimization programs with the same type of stochastic oracle.  Our abstraction unifies the analysis of actor-critic methods in reinforcement learning; we show how our main results can be specialized to recover the best-known convergence rate for policy optimization for an infinite-horizon average-reward Markov decision process (MDP) and to derive state-of-the-art rates for the online linear-quadratic regulator (LQR) controller and policy optimization using entropy regularization.

Our optimization framework is described as follows.  The overall goal is to solve  the (possibly non-convex)  program
\begin{align}
	\theta^{\star}=\argmin_{\theta\in\mathbb{R}^d}~f(\theta),
	\label{prob:obj}
\end{align}
where the gradient of $f$ is accessed through a stochastic oracle $H(\theta,\omega,X)$.  The three arguments to $H$ are the decision variable $\theta$, an auxiliary variable $\omega\in\mathbb{R}^r$, and a random variable $X$ drawn over a compact sample space $\Xcal$.  For a fixed $\theta$, there is a single setting  of the auxiliary variable, which we will denote $\omega^{\star}(\theta)$, such that $H$ returns an unbiased estimate of the gradient $\nabla f(\theta)$ when $X$ is drawn from a particular distribution $\mu_\theta$, 
\begin{align}
	\nabla f(\theta) = \Eset_{X\sim\mu_{\theta}}[H(\theta,\omega^{\star}(\theta),X)].
	\label{operator:H}
\end{align}
For other settings of $\omega\neq\omega^{\star}(\theta)$ or $X$ drawn from a distribution other than $\mu_\theta$, $H(\theta,\omega,X)$ will be (perhaps severely) biased.  The mapping from $\theta$ to the ``optimal'' setting of the auxiliary variable $\omega^{\star}(\theta)$ is implicit; it is determined by solving a nonlinear system of equations defined by another stochastic sampling operator $G:\mathbb{R}^d\times\mathbb{R}^r\times\Xcal\rightarrow\mathbb{R}^r$.  Given $\theta$, $\omega^{\star}(\theta)$ is the (unique, we assume) solution to
\begin{align}
	\Eset_{X\sim\mu_{\theta}}[G(\theta,\omega^{\star}(\theta),X)] = 0.
	\label{operator:G}
\end{align}
Combining \eqref{operator:H} and \eqref{operator:G}, solving \eqref{prob:obj} is equivalent to finding $(\theta^{\star},\omega^{\star}(\theta^{\star}))$ such that
\begin{align}\label{prob:GH=0}
	\left\{\begin{array}{ll}
		&\hspace{-0.5cm}\Eset_{X\sim\mu_{\theta^{\star}}}[H(\theta^{\star},\omega^{\star}(\theta^{\star}),X)] = 0,\\
		&\hspace{-0.5cm}\Eset_{X\sim\mu_{\theta^{\star}}}[G(\theta^{\star},\omega^{\star}(\theta^{\star}),X)] = 0.
	\end{array}\right.    
\end{align}

In the applications we are interested in, we only have indirect access to the distribution $\mu_\theta$.  Instead of parameterizing a distribution directly, $\theta$ parameterizes a set of probability transition kernels on $\Xcal\times\Xcal$ through a mapping $\Pcal:\Xcal\times\mathbb{R}^d\rightarrow\text{dist}(\Xcal)$.  Given $\theta$ and $X$, we will assume that we can generate a sample $X'\sim\Pcal(\cdot|X,\theta)$ using one of these kernels. Each of these $\Pcal(\cdot|\cdot,\theta)$ induces a different Markov chain and a different stationary distribution $\mu_\theta$, which is what is used in \eqref{operator:H} and \eqref{operator:G} above. 

This problem structure is motivated by online algorithms for reinforcement learning.
In this class of problems, which we describe in detail in Section~\ref{sec:applications}, we are searching for a control policy parameterized by $\theta$ that minimizes a long term cost captured by the function $f$.  The gradient for this cost depends on the value function under the policy indexed by $\theta$, which is specified implicitly through the Bellman equation, analogous to \eqref{operator:G} above.  These problems often also have a mechanism for generating samples, either through experiments or simulations, that makes only implicit use of the transition kernel $\Pcal$.

If we were able to generate independent and identically distributed (i.i.d.) samples from $\mu_{\theta}$ given any $\theta$ then a straightforward approach to solving \eqref{prob:GH=0} would be to use a nested-iterations approach. With $\theta$ fixed, we could sample $X\sim\mu_{\theta}$ repeatedly to (approximately) solve \eqref{operator:G} using stochastic approximation, then with $\theta$ and $\omega^{\star}(\theta)$ fixed, sample again and average over evaluations of $H$ to provide an (almost unbiased) estimate of the gradient in \eqref{operator:H}.  While there is an established methodology for analyzing this type of nested loop algorithm, variants of which have appeared in the control and reinforcement learning literature in \cite{yang2018finite,yang2019global,qiu2021finite}, they tend to be ``sample hungry'' in practice.  The method we study below involves only a single loop and is truly online, updating estimates of both the auxiliary and decision variables after each observed sample. This type of update is usually favored in practice as the implementation is simpler and empirically requires fewer samples  \cite{xu2020non,wu2020finite}.

An additional complication is that for a given $\theta$, actually generating samples from the induced stationary distribution $\mu_\theta$ is expensive and often impractical, requiring many draws using the kernels $\Pcal$ for a burn-in period and re-starting the Markov chain many times.  Our online update simply draws a single sample in every iteration $k$ according to $X_{k}\sim\Pcal(\cdot|X_{k-1},\theta_{k})$, where $\theta_{k}$ is the decision variable iterate whose updates will be discussed in Section~\ref{sec:two-time-scale-SGD}. This results in a time-varying Markov chain for $\{X_k\}$, which requires a delicate analytical treatment. Nevertheless, we are able to show that this non-stationarity settles down over time and to establish rates of convergence of the decision and auxiliary variable iterates.

Finally, we note that along with applications in reinforcement learning, problems that involve the types of coupled system of equations described above also appear in bi-level optimization  \cite{hong2020two} and composite optimization \cite{wang2017stochastic,chen2021solving}.

\subsection{Main Contribution} The main contributions of this paper are twofold.

First, we identify the common structure shared by various two-time-scale reinforcement learning algorithms and provide a clear mathematical abstraction of these algorithms and their objectives.
We develop a new two-time-scale SGD method for solving this class of objectives under the general stochastic oracle described above. This is an iterative algorithm that requires exactly one evaluation of $H$, one evaluation of $G$, and one query of the sampling oracle in each iteration. We then analyze the finite-time (and equivalently finite-sample) complexity of the proposed method when the objective function has different types of structure, namely, strong convexity, the Polyak-\L ojasiewicz (P\L) condition, and general non-convexity. For each type of function structure we provide an explicit upper bound on the convergence rate that uses carefully tuned decay rates of the step sizes.
The finite-time convergence rates are summarized in Table~\ref{table:result_summary}, where $k$ denotes the number of algorithm iterations. Since we draw exactly one sample in every iteration of the algorithm, the sample complexity has the same order as the time complexity. The key technical challenge is handling the coupling between the two iterates while addressing the impact of time-varying Markovian randomness. 
We analyze a coupled dynamical system that characterizes the convergence of the two variables by constructing a Lyapunov function tailored to the function structure of each objective.  This allows us to perform the analysis in a highly modular fashion, making it easily extendable to other types of function structure including general convexity and the restricted secant inequality \cite{karimi2016linear}.


Second, we demonstrate how our theoretical results apply to three policy optimization problems in reinforcement learning.
Although these algorithms can all be described on a high level as actor-critic methods, they operate under significantly different structure including continuous/discrete state spaces, the dynamical model of the state transition, and the reward function; nevertheless, our framework and analysis apply easily to all cases.
We start with the infinite-horizon average-reward MDP setting with finite state and action spaces, where our analysis recovers the state-of-the-art result for online actor-critic algorithms for policy optimization.  Next we describe how our general framework can also be specialized as an online natural actor-critic method for solving the classic linear-quadratic regulator (LQR) problem, providing the first known finite-sample complexity guarantee for an online algorithm for this problem. 
Specifically, we prove that after $k$ iterations, this method finds a control that is within $\Ocal(k^{-2/3})$ of the optimal solution in terms of the infinite-horizon average cost.  
Key to our analysis is exploiting the fact that the objective of the LQR problem obeys the P\L~condition.
Finally, we consider the policy optimization problem in a infinite-horizon discounted-reward MDP under entropy regularization. In this case, the objective function again obeys the P\L~condition, and our framework translates to an actor-critic algorithm that converges with rate $\Ocal(k^{-2/3})$ to the globally optimal objective function. This algorithm and its analysis are both previously unknown in the literature.
\vspace{-10pt}

\begin{table}[!h]
\centering
\caption{\small\sl Summary of Main Results - Time and Sample Complexity.
}
\vspace{-.1cm}
\setlength{\tabcolsep}{5pt}
\begin{tabular}{cccccc}
        \toprule
        \makecell{Structural \\ Property of $f$}  &  Metric  &  Rate & \makecell{Order of\\ $\alpha_k$, $\beta_k$} & \makecell{Standard\\ SGD Rate} & \makecell{Applications\\ Presented in\\ Section}\\
        \midrule
        Strong Convexity & $\|\theta_k-\theta^{\star}\|^2$ & $\widetilde{\Ocal}(k^{-\frac{2}{3}})$ & $k^{-1}, k^{-\frac{2}{3}}$ & $\Ocal(k^{-1})$ & \ref{sec:TDC} \\

        P\L~Condition & $f(\theta_k)-f(\theta^{\star})$ & $\widetilde{\Ocal}(k^{-\frac{2}{3}})$ & $k^{-1}, k^{-\frac{2}{3}}$ & $\Ocal(k^{-1})$ & \ref{sec:applications_lqr}, \ref{sec:ac_entropy} \\
        Non-convexity & $\|\nabla f(\theta_k)\|^2$\hspace{-2pt} & $\widetilde{\Ocal}(k^{-\frac{2}{5}})$ & $k^{-\frac{3}{5}}, k^{-\frac{2}{5}}$ & $\Ocal(k^{-\frac{1}{2}})$ & \ref{sec:MDP_AC}\\
        \bottomrule
        \bottomrule
        \end{tabular}
\label{table:result_summary}
\end{table}

\subsection{Related Works}\label{sec:related_works}

Our work is closely related to the existing literature on two-time-scale stochastic approximation (SA), bi-level and composite optimization, actor-critic algorithms in reinforcement learning, and single-time-scale stochastic optimization algorithms under unbiased or biased (sub)gradients. In this section, we discuss the recent advances in these domains to give context to our contributions.

\vspace{.1in}

\noindent\textbf{Two-Time-Scale SA.} Two-time-scale SA solves a system of equations similar in form to \eqref{prob:GH=0}, but typically considers the setting where $\mu_{\theta}=\mu$ is independent of the decision variable $\theta$. The convergence of two-time-scale SA is traditionally established by analyzing an associated ordinary differential equation \cite{borkar2000ode}. Finite-time convergence of two-time-scale SA has been studied in the case where $H$ and $G$ are linear \cite{KondaT2004, DalalTSM2018, Dalal_Szorenyi_Thoppe_2020, DoanR2019,GuptaSY2019_twoscale,Kaledin_two_time_SA2020} and in more general nonlinear settings \cite{MokkademP2006,doan2021finite}, under either i.i.d.\ or Markovian samples. In these previous works, the analysis for the nonlinear setting is restricted to the case where $H$ and $G$ are both strongly monotone, while our work studies a wide range of function structures including strong convexity, the P\L~condition, and general non-convexity. \looseness=-1

\vspace{.1in}

\noindent\textbf{Bi-Level and Composite Optimization.} The optimization objective in our work is closely connected to the bi-level optimization framework \cite{colson2007overview,hong2020two,chen2022single} which solves programs structured as
\begin{align}
    \min_{x} f_1(x,y^{\star}(x)) \quad\text{subject to  }y^{\star}(x)\in\argmin_y f_2(x,y).\label{eq:bilevel_opt}
\end{align}
From the first-order optimality condition, \eqref{eq:bilevel_opt} is equivalent to finding a stationary point $(x',y')$ that observes
\begin{align*}
    \nabla_x f_1(x',y')=0,\quad \nabla_y f_2(x',y')=0.
\end{align*}
This is a special case of our objective in \eqref{prob:GH=0} where $G$ is a gradient mapping. However, in RL applications $G$ usually abstracts the Bellman backup operator which is associated with the estimation of the value function. It is well-known that the Bellman backup operator is not the gradient of any function. In this sense, our framework is more general and suitable for modeling algorithms in RL. In addition, the bi-level optimization literature considers a stochastic oracle with a fixed distribution $\mu$, while we solve \eqref{prob:GH=0} with the distribution of the samples also depending on the decision variable. This is another important generalization as many realistic problems in control and reinforcement learning can only be abstracted as \eqref{prob:GH=0} with $\mu_{\theta}$ being a function of the control variable $\theta$. Making this generalization requires generating decision-variable-dependent samples from a Markov chain whose stationary distribution shifts over iterations, which creates additional challenges in the analysis.\looseness=-1

{\color{blue}
We also note the connection of our objective to stochastic composite optimization \cite{wang2017stochastic,chen2021solving}, which solves optimization problems of the form
\begin{align}
    \min_{x}g_1(g_2(x)).\label{eq:composite_opt}
\end{align}
At a first glance, \eqref{eq:composite_opt} reduces to \eqref{eq:bilevel_opt} by choosing $g_1=f_1$ and $g_2(x)=(x,y^{\star}(x))$ where $y^{\star}(x)$ is the minimizer of $f_2(x,\cdot)$. 
However, the key assumption made in these works, and in the stochastic composite optimization literature in general, is the differentiability of $g_2$ (and $g_1$) and access to an oracle that returns the stochastic gradient $\nabla g_2$. This assumption is restrictive in reinforcement learning where only indirect information about $g_2$ is available.
}

\vspace{.1in}

\noindent\textbf{Actor-Critic Algorithms.} In the RL literature the aim of actor-critic algorithms is also to solve a problem similar to \eqref{prob:GH=0}, where $\theta$ and $\omega^{\star}(\theta)$ are referred to as the actor and critic, respectively; {\color{blue} see for example \cite{qiu2021finite, kumar2019sample,xu2020non, wu2020finite}, among which \cite{qiu2021finite, kumar2019sample,xu2020non} all employ nested-loop updates and require constant state resetting, making these algorithms inconvenient to implement in practice. The authors in \cite{wu2020finite} consider an online, single-loop setting similar to the one studied in our paper. 
In fact, the algorithm in \cite{wu2020finite} is a special case of our framework with a non-convex objective function. Our analysis recovers the result of \cite{wu2020finite} while slightly loosening the assumptions --- by carefully exploiting a "strong monotonicity" property of the operator $G$ we are able to remove a projection operator used by \cite{wu2020finite} to limit the growth of the critic parameter.
The same type of projection operator is also used in \cite{xu2020non}, which otherwise makes similar assumptions to ours including the Lipschitz continuity/smoothness of the objective function and the uniform ergodicity of the Markov chain (Assumption~\ref{assump:markov-chain}).}

\vspace{.1in}

\noindent\textbf{Single-Variable Stochastic Optimization.} When the samples are i.i.d., stochastic gradient/subgradient algorithms are fairly well-understood for smooth (see \cite{amari1993backpropagation,gower2019sgd,khaled2020tighter} and the references therein) and non-smooth \cite{ruszczynski1987linearization,boyd2008stochastic,ruszczynski2020convergence} functions.
In the smooth setting, \cite{chen2019finite,sun2018markov,zeng2023finite} study various SGD/SA algorithms under samples generated from time-invariant state transition probabilities (we will later refer to this as a time-invariant Markov chain) and show that the convergence rates are only different from that under i.i.d.\ samples by a logarithmic factor. The key argument used in these works is that the Markovian samples behave similarly to i.i.d.\ samples on a mildly dilated time scale.

In many policy optimization algorithms in RL, the samples are drawn under the control of the current policy. As the policy gets updated, the transition probabilities shift, resulting in a Markov chain with a time-varying stationary distribution (we will refer to this as the time-varying Markov chain). This setting requires more sophisticated mathematical treatment.
The single-variable SA algorithm under time-varying Markovian samples is first analyzed by \cite{zou2019finite}, while our paper is among the first works to extend the analysis to the scenario where two coupled variables are updated simultaneously. An important assumption to treat time-varying Markovian samples is the ``Lipschitz continuity'' of the transition probability with respect to the policy, which we note holds in Markov decision processes and linear-quadratic regulator discussed in Section~\ref{sec:applications}.

\section{Two-Time-Scale Stochastic Gradient Descent}\label{sec:two-time-scale-SGD}
In this section, we present our two-time-scale SGD method (formally stated in Algorithm \ref{Alg:two-time-scale-SGD}) for solving \eqref{prob:GH=0} under the gradient oracle discussed in Section \ref{sec:introduction}. 
In the algorithm, $\theta_{k}$ and $\omega_{k}$ are estimates of $\theta^{\star}$ and $\omega^{\star}(\theta^{\star})$. The random variables $\{X_{k}\}$ are generated by a Markov process parameterized by $\{\theta_{k}\}$ under the transition kernel $\Pcal$, i.e.,   
\begin{align}
    X_0 \stackrel{\theta_1}{\longrightarrow} X_1 \stackrel{\theta_2}{\longrightarrow} X_2 \stackrel{\theta_3}{\longrightarrow} \cdots \stackrel{\theta_{k-1}}{\longrightarrow} X_{k-1} \stackrel{\theta_{k}}{\longrightarrow} X_{k}.
    \label{eq:Markovchain_timevarying}
\end{align}
As $\theta_k$ changes in every iteration, so do the dynamics of this Markov process that generates the data.  At a finite step $k$, $X_k$ is in general not an i.i.d.\ sample from the stationary distribution $\mu_{\theta_k}$, implying that $H(\theta_{k},\omega_{k},X_k)$ employed in the update \eqref{alg:update_theta} is not an unbiased estimate of $\nabla f(\theta_k)$ even if $\omega_k$ tracks $\omega^{\star}(\theta_{k})$ perfectly. This sample bias, along with the gap between $\omega_{k}$ and $\omega^{\star}(\theta_{k})$, affects the variables $\theta_{k+1}$ and $\omega_{k+1}$ of the next iteration and accumulates inaccuracy over time which needs a careful treatment.\looseness=-1

\begin{algorithm}
\caption{Two-Time-Scale Stochastic Gradient Descent}
\label{Alg:two-time-scale-SGD}
\begin{algorithmic}[1]
\STATE{\textbf{Initialization:} decision variable $\theta_0$, auxiliary variable $\omega_0$, sample $X_0$, step size sequences $\{\alpha_k\}$ for decision variable update, $\{\beta_k\}$ for auxiliary variable update}
\FOR{$k=1,2,3,...$}
\STATE{Decision variable update:
        \vspace{-5pt}
        \begin{align}
        \theta_{k+1} &= \theta_{k} - \alpha_k H(\theta_k,\omega_k,X_{k})
        \label{alg:update_theta}
        \end{align}}
\STATE{\vspace{-15pt} Auxiliary variable update:
        \vspace{-5pt}
        \begin{align}
        &\omega_{k+1}=\omega_k-\beta_k G(\theta_{k+1},\omega_k, X_k)
        \label{alg:update_omega}
        \end{align}
        }
\STATE{\vspace{-15pt} Draw sample 
\vspace{-10pt}
\begin{align}
    X_{k+1}\sim \Pcal(\cdot\mid X_{k}, \theta_{k+1})
    \label{eq:update_sample}
\end{align}}
\ENDFOR
\end{algorithmic}
\end{algorithm}

The updates use two different step sizes, $\alpha_{k}$ and $\beta_{k}$. We choose $\alpha_{k}\ll\beta_{k}$ as a way to approximate, very roughly, the nested-loop algorithm that runs multiple auxiliary variable updates for each decision variable update.  Many small critic updates get replace with a single large one.  In other words, the auxiliary variable $\omega_{k}$ is updated at a faster time scale (larger step size) as compared to $\theta_{k}$ (smaller step size). 

The ratio $\beta_{k}/\alpha_{k}$ can be interpreted as the time-scale difference.  We will see that this ratio needs to be carefully selected  based on the structural properties of the function $f$ for the algorithm to achieve the best possible convergence.  Table~\ref{table:result_summary} provides a brief summary of our main theoretical results, which characterizes the finite-time complexity of Algorithm \ref{Alg:two-time-scale-SGD} and the corresponding optimal choice of step sizes under different function structures. Table~\ref{table:result_summary} also contrasts the convergence of Algorithm~\ref{Alg:two-time-scale-SGD} with the rates of standard SGD (which in our context means that the samples are i.i.d.\ and the auxiliary variable is always exactly accurate).

\section{Motivating Applications}\label{sec:applications}

In this section, we show how our results on two-time-scale optimization apply to a variety of policy evaluation and optimization algorithms in RL. The first three applications can be categorized as actor-critic algorithms for policy optimization.
The objectives are non-convex in these applications, but the second and third problems are more structured and observe the P\L~condition. In Section~\ref{sec:TDC}, we briefly discuss an application of the framework to two-time-scale gradient-based policy evaluation algorithms where the objective function is strongly convex.\looseness=-1

\subsection{Online Actor-Critic Method for MDPs (Non-Convex)}\label{sec:MDP_AC}
We consider the standard infinite-horizon average-reward MDP model $\Mcal=(\Scal,\Acal,\Pcal,r)$, where $\Scal$ is the state space, $\Acal$ is the action space, $\Pcal:\Scal\times\Acal\rightarrow\Delta_{\Scal}$ denotes the transition probabilities, and $r:\Scal\times\Acal\rightarrow[-1,1]$ is the reward. Our aim is to find the policy $\pi_{\theta}\in\Delta_{\Acal}^{\Scal}$, parameterized by $\theta\in\mathbb{R}^d$ (where $d$ may be much smaller than $|\Scal|\hspace{-2.5pt}\times\hspace{-2.5pt}|\Acal|$), that maximizes the average cumulative reward\vspace{-5pt}
\begin{align}
    \theta^{\star}=\argmax_{\theta\in\mathbb{R}^d} J(\theta)\triangleq\lim_{K \rightarrow \infty}\frac{1}{K}\Eset\Big[\sum_{k=0}^{K} r\left(s_{k}, a_{k}\right)\Big]=\mathbb{E}_{s \sim \mu_{\theta}, a \sim \pi_{\theta}}[r(s, a)],
    \label{eq:actorcritic_RL_obj}
\end{align}
where $\mu_{\theta}$ denotes the stationary distribution of the states induced by the policy $\pi_{\theta}$. Defining the (differential) value function of the policy $\pi_{\theta}$\vspace{-5pt}
\begin{align*}
    V^{\pi_{\theta}}(\cdot)=\mathbb{E}\Big[\sum_{k=0}^{\infty}\left(r\left(s_{k}, a_{k}\right)-J(\theta)\right) \mid s_{0}=\cdot\Big],
\end{align*}
we can use the well-known policy gradient theorem 
to express the gradient of the objective function in \eqref{eq:actorcritic_RL_obj} as\vspace{-5pt}
\begin{align*}
    \nabla J(\theta)\hspace{-2pt}=\hspace{-2pt}\mathbb{E}_{s \sim \mu_{\theta}(\cdot), a \sim \pi_{\theta}(\cdot \mid s), s'\sim\Pcal(\cdot\mid s,a)}\Big[\hspace{-3pt}\left(r(s,a)\hspace{-2pt}-\hspace{-2pt}J(\theta)\hspace{-2pt}+\hspace{-2pt}V^{\pi_{\theta}}(s')\hspace{-2pt}-\hspace{-2pt}V^{\pi_{\theta}}(s)\right)\nabla \hspace{-2pt}\log \pi_{\theta}(a \mid s)\Big].
\end{align*}

Optimizing \eqref{eq:actorcritic_RL_obj} with (stochastic) gradient ascent methods requires evaluating $V^{\pi_{\theta}}$ and $J(\theta)$ at the current iterate of $\theta$, which are usually unknown and/or expensive to compute exactly. 
``Actor-critic'' algorithms attack this problem on two scales as discussed in the sections above: an actor keeps a running estimate of the policy parameters $\theta_k$ ,
while a critic approximately tracks the differential value function for $\theta_k$ to aid the evaluation of the policy gradient. 

For problems with large state spaces, it is often necessary to use a low-dimensional parameter $\omega\in\mathbb{R}^m$ to approximate $V^{\pi_{\theta}}$ where $m\ll |\Scal|$.
In this work, we consider the linear function approximation setting where each state $s$ is encoded by a feature vector $\phi(s)\in\mathbb{R}^m$ and the approximate value function is $\widehat{V}^{\pi_{\theta},\psi}(s)=\phi(s)^{\top}\psi$. Under the assumptions that the Markov chain of the states induced by any policy is uniformly ergodic (equivalent of Assumption \ref{assump:markov-chain} in this context) and that the feature vectors $\{\phi(s)\}_{s\in\Scal}$ are linearly independent, it can be shown that a unique optimal pair $(J(\theta),\psi^{\star}(\theta))$ exists that solves the projected Bellman equation\vspace{-3pt}
\begin{align*}
    \mathbb{E}_{s \sim \mu_{\theta}(\cdot), a \sim \pi_{\theta}(\cdot \mid s), s'\sim\Pcal(\cdot\mid s,a)}\left[\hspace{-3pt}\begin{array}{c}
    J(\theta)-r(s,a)\\
    \left(r(s,a)-J(\theta)+\phi(s')^{\top}\psi^{\star}(\theta)-\phi(s)^{\top}\psi^{\star}(\theta)\right)\phi(s)
    \end{array}\hspace{-3pt}\right]=0.
\end{align*}
We use an auxiliary variable $\omega=(\hat{J},\psi)$ to track the solution to this Bellman equation.

Due to the limit in the representational power of the function approximation, there is an approximation error between $V^{\pi_{\theta}}$ and $\hat{V}^{\pi_{\theta},\psi^{\star}(\theta)}$ as a function of $\theta$, which we define over the stationary distribution\vspace{-3pt}
\begin{align*}
    \epsilon_{\text{approx}}(\theta)=\sqrt{\mathbb{E}_{s \sim \mu_{\theta}}\left[(\phi(s)^{\top} \psi^{\star}(\theta)-V^{\pi_{\theta}}(s))^{2}\right]}.
\end{align*}
We assume the existence of a constant $\epsilon_{\text{approx}}^{\max}$ such that $\epsilon_{\text{approx}}(\theta)\hspace{-2pt}\leq\hspace{-2pt}\epsilon_{\text{approx}}^{\max}$ for all $\theta\in\mathbb{R}^d$.\looseness=-1

Comparing this problem with Eq. \eqref{operator:H} and \eqref{operator:G}, it is clear that this is a special case of our optimization framework with $X=(s,a,s')$, $\omega^{\star}(\theta)=(J(\theta),\psi^{\star}(\theta))$ and
\begin{align*}
    &f(\theta) = -J(\theta),\,\,\, G(\theta,\omega,X) = [ \hat{J}-r(s,a), (r(s,a)-\hat{J}+\phi(s')^{\top}\psi-\phi(s)^{\top}\psi)\phi(s)^{\top}]^{\top},\notag\\
    &H(\theta,\omega,X) = -(r(s,a)-\hat{J}+\phi(s')^{\top}\psi-\phi(s)^{\top}\psi+\varepsilon_{\text{approx}}(\theta))\nabla \log \pi_{\theta}(a \mid s),
\end{align*}
where $\varepsilon_{\text{approx}}(\theta)$ is an error in the gradient of the actor carried over from the approximation error of the critic which can be upper bounded by $2\epsilon_{\text{approx}}^{\max}$ in expectation.
In this case, the function $-J(\theta)$ is non-convex and our two-time-scale SGD algorithm is guaranteed to find a stationary point of the objective function with rate $\widetilde{\Ocal}(k^{-2/5})$, up to errors proportional to $\epsilon_{\text{approx}}^{\max}$. This rate matches the state-of-the-art result derived in \cite{wu2020finite}. 
A subtle improvement of our analysis is that we do not need to perform the projection of the critic parameter onto a compact set that \cite{wu2020finite} requires in every iteration of the algorithm to guarantee the stability of the critic.



\subsection{Online Natural Actor-Critic Algorithm for LQR (P\L~Condition)}
\label{sec:applications_lqr}
\sloppy
In this section, we consider the infinite-horizon average-cost LQR problem
\begin{align}
\begin{aligned}
&\underset{\{u_{k}\}}{\text{minimize}}\quad \lim_{T\rightarrow\infty}\frac{1}{T}\mathbb{E}\Big[\sum_{k=0}^{T}\left(x_{k}^{\top} Q x_{k}+u_{k}^{\top} R u_{k}\right)\mid x_{0}\Big] \\
&\text{subject to}\quad  x_{k+1}=A x_{k}+B u_{k}+ w_k,
\end{aligned}
\label{eq:obj_LQR}
\end{align}
where $x_k\in\mathbb{R}^{d_1}$, $u_k\in\mathbb{R}^{d_2}$ are the state and the control variables, $w_k\sim N(0,\Psi)\in\mathbb{R}^{d_1}$ is time-invariant system noise, $A\in\mathbb{R}^{d_1\times d_1}$ and  $B\in\mathbb{R}^{d_1\times d_2}$ are the system transition matrices, and $Q\in\mathbb{R}^{d_1\times d_1}, R\in\mathbb{R}^{d_2\times d_2}$ are positive-definite cost matrices. 
It is well-known (see, for example, \cite[Chap.\ 3.1]{bertsekas2012dynamic}) that the optimal control sequence $\{u_k\}$ that solves \eqref{eq:obj_LQR} is a time-invariant linear function of the state 
\begin{align}
    u_k^\star = -K^{\star} x_k,\label{eq:lqr_u*}
\end{align}
where $K^{\star}\in\mathbb{R}^{d_2 \times d_1}$ is a matrix that depends on the problem parameters $A,B,Q,R$. This fact will allow us to reformulate the LQR as an optimization program over the feedback gain matrix $K$.  It is also true that optimizing over the set of stochastic controllers\looseness=-1
\begin{align*}
u_{k} = -K x_{k} + \sigma \epsilon_{k}, \quad \epsilon_{k} \sim \text{N}(0,\sigma^2\mathbf{I}),
\end{align*}
with $\sigma \geq 0$ fixed will in the end yield the same optimal $K^{\star}$ \cite{gravell2020learning}.  In the RL setting considered below, we will optimize over this class of stochastic controller as it encourages exploration.  Defining $\Psi_{\sigma}=\Psi+\sigma^2 BB^{\top}$, we can re-express the LQR problem as
\begin{align}
\begin{aligned}
    \underset{K}{\text{minimize }} & J(K)\triangleq \operatorname{trace}(P_{K}\Psi_{\sigma})+\sigma^2\operatorname{trace}(R)\\
    \text{s.t. }&P_{K}=Q+K^{\top} R K+(A-B K)^{\top} P_{K}(A-B K).
\end{aligned}
\label{eq:obj_K_LQR}
\end{align}

Our goal is to solve \eqref{eq:obj_K_LQR} when the system transition matrices $A$ and $B$ are unknown\footnote{We do assume, however, that we know the cost matrices $Q$ and $R$ or at least we can compute $x^{\top}Qx+u^{\top}Ru$ for any $x$ and $u$.} and we take online samples from a single trajectory of states $\{x_k\}$ and control inputs $\{u_k\}$. This problem has been considered recently in \cite{yang2019global}, and in fact much of our formulation is modeled on this work. The essential difference is that while \cite{yang2019global} works in the ``batch'' setting, where multiple trajectories are drawn for a fixed feedback gain estimate, our algorithm is entirely online.

Given a feedback gain $K$, we define
\[E_K\triangleq2\left(R+B^{\top} P_{K} B\right) K-2B^{\top} P_{K} A.\]
It turns out that the natural policy gradient of the objective function in \eqref{eq:obj_K_LQR}, which we denote by $\widetilde{\nabla} J$, is $\widetilde{\nabla} J(K) = E_K$.

To track $E_K$ when $A$ and $B$ are unknown it suffices to estimate $R+B^{\top} P_K B$ and $B^{\top} P_K A$. We define \vspace{-3pt}
\begin{align}
    \Omega_{K}=\left(\begin{array}{cc}
    \Omega_{K}^{11} & \Omega_{K}^{12} \\
    \Omega_{K}^{21} & \Omega_{K}^{22}
    \end{array}\right)=\left(\begin{array}{cc}
    Q+A^{\top} P_{K} A & A^{\top} P_{K} B \\
    B^{\top} P_{K} A & R+B^{\top} P_{K} B
    \end{array}\right),
\end{align}
of which $R+B^{\top} P_K B$ and $B^{\top} P_K A$ are sub-matrices.
We define the operator $\text{svec}(\cdot)$ as the vectorization of the upper triangular sub-matrix of a symmetric matrix with off-diagonal entries weighted by $\sqrt{2}$, and define $\text{smat}(\cdot)$ as the inverse of $\text{svec}(\cdot)$. We also define $\phi(x, u)=\text{svec}(\left[x^{\top},u^{\top}\right]^{\top}\left[x^{\top},u^{\top}\right])$ for any $x\in\mathbb{R}^{d_1},u\in\mathbb{R}^{d_2}$.
Then, it can be shown that $\Omega_{K}$ and $J(K)$ jointly satisfy the Bellman equation\vspace{-3pt}
\begin{align}
\hspace{-0.3cm}\mathbb{E}_{x\sim\mu_K,u\sim N(-K x,\sigma^2 I)}\left[M_{x,u,x',u'}\right]\left[\begin{array}{c}
J(K) \\
\text{svec}(\Omega_K)
\end{array}\right]=\mathbb{E}_{x\sim\mu_K,u\sim N(-K x,\sigma^2 I)}\left[c_{x,u}\right],
\label{eq:lqr_Goperator}
\end{align}
where the matrix $M_{x,u,x',u'}$ and vector $c_{x,u}$ are\vspace{-3pt}
\begin{gather*}
M_{x,u,x',u'}\hspace{-2pt}=\hspace{-2pt}\left[\hspace{-4pt}\begin{array}{cc}
1 & 0 \\
\phi(x, u) & \hspace{-4pt}\phi(x, u)\left[\phi(x, u)-\phi\left(x', u'\right)\right]^{\top}
\end{array}\hspace{-6pt}\right]\hspace{-2pt},\,\,c_{x,u}\hspace{-2pt}=\hspace{-2pt}\left[\hspace{-4pt}\begin{array}{c}
x^{\top} Q x+u^{\top} R u \\
\left(x^{\top} Q x+u^{\top} R u\right)\hspace{-2pt}\phi(x, u)
\end{array}\hspace{-7pt}\right]\hspace{-2pt}.
\end{gather*}
The solution to \eqref{eq:lqr_Goperator} is unique if $K$ is stable with respect to $A$ and $B$ \cite{yang2019global}. An auxiliary variable $\hat{\Omega}$ can be introduced to track $\Omega_K$ for the decision variable $K$.

We connect this to our optimization framework by noting that Eq.~\eqref{eq:lqr_Goperator} corresponds to Eq.~\eqref{operator:G}  with $K$, $(x,u,x',u')$, and $[J(K),\text{svec}(\Omega_K)^{\top}]^{\top}$ mirroring $\theta$, $X$, and $\omega^{\star}(\theta)$, respectively. The natural gradient oracle in this case is $H(\theta,\omega,X)=2\hat{\Omega}^{22}K-2\hat{\Omega}^{21}$, which does not depend on the samples $X$ directly, and the operator $G$ is $G(\theta,\omega,X)=-M_{x,u,x',u'}\omega+c_{x,u}$. A key structure of \eqref{eq:obj_LQR} is that the objective function is non-convex but observes the P\L~condition \cite{yang2019global}, which we formally discuss later in Assumption~\ref{assump:PL_condition}. As a result, applying Algorithm~\ref{Alg:two-time-scale-SGD} to this problem leads to an online actor-critic flavored algorithm that converges with rate $\widetilde{\Ocal}(k^{-2/3})$ under proper assumptions.
To our best knowledge, our work is the first to study the online actor-critic method for solving the LQR, and our result vastly improves over the rate $\widetilde{\Ocal}(k^{-1/5})$ of the nested-loop actor-critic algorithm derived in \cite{yang2019global} which also operates under more restrictive assumptions (e.g.\ sampling from the stationary distribution, boundedness of the iterates). 

Section~\ref{sec:lqr_details} presents a more detailed description of the online actor-critic algorithm, gives a formal statement of the assumptions, states the theoretical results carefully, and gives a small-scale numerical simulation. After the appearance of our paper, a recent work \cite{zhou2022single} improves the convergence rate further to $\widetilde{\Ocal}(k^{-1})$ by introducing a single-time-scale actor-critic algorithm where the critic carries out least square temporal difference learning \cite{bradtke1996linear}. However, we note that their algorithm generates multiple trajectories to estimate the gradient in each critic update and thus requires constant resetting of the agent to desired states, which may not be feasible in many practical applications.

\subsection{Online Actor-Critic Method for Regularized MDPs (P\L~Condition)}\label{sec:ac_entropy}

As a third application of our framework,
we study the policy optimization problem for the infinite-horizon discounted-reward MDP $\Mcal=(\Scal,\Acal,\Pcal,r,\gamma)$ where $\gamma\in(0,1)$ is the discount factor and the rest are defined in the same manner as above in Section~\ref{sec:MDP_AC}.
We restrict our attention to the tabular setting where the parameter $\theta$ encodes the policy through the softmax function
\begin{align*}
    \pi_{\theta}(a \mid s)&=\frac{\exp \left(\theta_{s, a}\right)}{\sum_{a' \in \Acal} \exp \left(\theta_{s, a'}\right)}.
 \end{align*}
 
To accelerate the convergence of the actor-critic algorithm, we regularize the objective with the policy entropy as proposed by \cite{mei2020global}. Specifically, with regularization weight $\tau>0$, the regularized value function of a policy $\pi$ is
\begin{align*}
    V_{\tau}^{\pi}(s)&=\mathbb{E}_{a_k\sim\pi(\cdot\mid s_k),s_{k+1}\sim\Pcal(\cdot\mid s_k,a_k)}\Big[\sum_{k=0}^{\infty}\gamma^k\big(r(s_k,a_k)-\tau\log\pi(a_k\mid s_k)\big)\mid s_0=s\Big].
\end{align*}
Under the initial state distribution $\rho\in\Delta_{\Scal}$, the expected cumulative reward collected by policy $\pi$ is $J_{\tau}(\pi)=\mathbb{E}_{s\sim\rho}[V_{\tau}^{\pi}(s)]$.
We consider solving the policy optimization problem\looseness=-1 
\begin{align*}
    \max_{\pi}J_{\tau}(\pi).
\end{align*}
Expressing the gradient of the objective with the policy gradient theorem, we have
\begin{align*}
    &\nabla_{\theta} J_{\tau}(\pi_{\theta})=\frac{1}{1-\gamma}\mathbb{E}\big[\left(r(s,a)-\tau\log\pi_{\theta}(a\mid s)+\gamma V_{\tau}^{\pi_{\theta}}(s')-V_{\tau}^{\pi_{\theta}}(s)\right)\nabla_{\theta}\log\pi_{\theta}(a\mid s)\big],
\end{align*}
where the expectation is taken over $s\sim d_{\rho}^{\pi_{\theta}},a\sim\pi_{\theta}(\cdot\mid s),s'\sim\Pcal(\cdot\mid s,a)$, and the the discounted visitation distribution $d_{\rho}^{\pi}\in\Delta_{\Scal}$ is defined as
\begin{align*}
    d_{\rho}^{\pi}(s)=(1-\gamma)\mathbb{E}_{a_k\sim\pi(\cdot\mid s_k),s_{k+1}\sim\Pcal(\cdot\mid s_k,a_k)}[\sum\nolimits_{k=0}^{\infty} \gamma^k \1(s_k=s) \mid s_0\sim\rho].
\end{align*}

To evaluate the gradient, we need to compute the regularized value function $V_{\tau}^{\pi_{\theta}}$, which is the solution to the following Bellman equation
\begin{align*}
    \mathbb{E}_{s\sim d_{\rho}^{\pi_{\theta}},a\sim\pi_{\theta}(\cdot\mid s),s'\sim\Pcal(\cdot\mid s,a)}\Big[r(s,a)-\tau\log\pi_{\theta}(a\mid s)+\gamma V_{\tau}^{\pi_{\theta}}(s')-V_{\tau}^{\pi_{\theta}}(s)\Big]=0.
\end{align*}
Interestingly, \cite{kondathesis} shows that we can regard $d_{\rho}^{\pi}$ as the stationary distribution under $\pi$ in an environment with the modified transition probability
\begin{align*}
    \widetilde{P}(\cdot\mid s,a)=\gamma P(\cdot\mid s,a)+(1-\gamma)\rho(\cdot).
\end{align*}
This observation allows us to generate Markovian samples $(s,a,s')$ with $d_{\rho}^{\pi_{\theta}}\otimes \pi_{\theta}\otimes \Pcal$ as the stationary distribution, in an online manner that resembles \cite{barakat2022analysis}[Algorithm 1].\looseness=-1

In the actor-critic framework, we introduce a critic (auxiliary variable) $\widehat{V}$ to estimate the solution of the Bellman equation under the current policy iterate. Our optimization framework abstracts this problem by choosing
\begin{align*}
    &X=(s,a,s'),\,\,\,\omega=\widehat{V},\,\,\,f(\theta) = -J_{\tau}(\pi_{\theta}),\notag\\
    &H(\theta,\omega,X) = \frac{1}{1-\gamma}(r(s,a)-\tau\log\pi_{\theta}(a\mid s)+\gamma\widehat{V}(s')-\widehat{V}(s))\nabla\log\pi_{\theta}(a\mid s),\notag\\
    &G(\theta,\omega,X) =  r(s,a)-\tau\log\pi_{\theta}(a\mid s)+\gamma\widehat{V}(s')-\widehat{V}(s).
\end{align*}
The objective function is non-convex but satisfies the P\L~condition under standard assumptions (see \cite{mei2020global}[Lemma 15]). Our two-time-scale SGD framework specializes to an online actor-critic algorithm, which by our analysis to be discussed in Section~\ref{sec:conv_pl} is guaranteed to find a globally optimal solution of the regularized objective with rate $\widetilde{\Ocal}(k^{-2/3})$.
To our best knowledge, this is the first time such data-driven algorithms are studied for solving an entropy-regularized MDP in the tabular setting.
Compared with the result presented in Section~\ref{sec:MDP_AC}, the introduction of the entropy regularization leads to an accelerated convergence rate. We note that the gap between the solutions to the regularized and original MDP is proportional to the regularization weight $\tau$ \cite{cen2021fast,zeng2022regularized}. By carefully choosing $\tau$, solving the regularized MDP provides a reliable and efficient way to find the approximate solution of the original unregularized MDP.

\subsection{Two-Time-Scale Policy Evaluation Algorithms (Strongly Convex)}\label{sec:TDC}
Our framework also abstracts GTD (gradient
temporal-difference), GTD2
, and TDC (temporal difference learning with gradient correction)
algorithms \cite{sutton2008convergent,sutton2009fast}, which are gradient-based two-time-scale policy evaluation algorithms in RL. They can be viewed as degenerate special cases of our framework where the expectation in \eqref{prob:GH=0} is taken over a fixed distribution $\mu$ that does not depend on $\theta$, and therefore do not require the full capacity of our framework. The objective function in this problem is strongly convex, and our framework under proper assumptions guarantees a convergence rate of $\widetilde{\Ocal}(k^{-2/3})$, which matches the analysis in \cite{xu2019two}. As this subject is well-studied, we skip the detailed discussion of the problem formulation and algorithm statement and refer interested readers to \cite{sutton2008convergent,sutton2009fast,xu2019two}.

\section{Technical Assumptions}
\label{sec:assumptions}
In this section, we present the main technical assumptions important in our later analysis. We first consider the Lipschitz continuity of $H$ and $G$. 
\begin{assump}\label{assump:HG_smooth}
There exists a constant $L > 0$ such that for all $\theta_{1},\theta_{2}\in\Rset^{d},$ $\omega_{1},\omega_{2}\in\Rset^{r}$, and $X\in\Xcal$ 
\begin{align}
    &\left\|H(\theta_1,\omega_1,X)-H(\theta_2,\omega_2,X)\right\| \leq L\left(\|\theta_1-\theta_2\|+\|\omega_1-\omega_2\|\right),\notag\\
    &\left\|G(\theta_1,\omega_1,X)-G(\theta_2,\omega_2,X)\right\| \leq L\left(\|\theta_1-\theta_2\|+\|\omega_1-\omega_2\|\right).\label{assump:HG_smooth:ineq}
\end{align}
\end{assump}
We also assume that the objective function $f$ is smooth. 
\begin{assump}\label{assump:f_smooth}
There exists a constant $L>0$ such that for all $\theta_{1},\theta_{2}\in\Rset^{d}$
\begin{align}
    \left\|\nabla f(\theta_1)-\nabla f(\theta_2)\right\| \leq L\|\theta_1-\theta_2\|.\label{assump:f_smooth:ineq}
\end{align}
\end{assump}
Assumptions \ref{assump:HG_smooth} and \ref{assump:f_smooth} are common in the literature of stochastic approximation \cite{karmakar2018two,doan2021finite} and hold in the actor-critic methods discussed in Section \ref{sec:applications}. Next, we assume that the operator $G(\theta,\cdot,X)$ is strongly monotone in expectation at $\omega^{\star}(\theta)$ (which we have assumed is unique).

\begin{assump}\label{assump:stronglymonotone_G}
There exists a constant $\lambda>0$ such that
\begin{align*}
    \left\langle \Eset_{X\sim\mu_{\theta}}[G(\theta,\omega,X)],\omega-\omega^{\star}(\theta)\right\rangle \leq -\lambda\|\omega-\omega^{\star}(\theta)\|^2,\quad \forall \theta \in \mathbb{R}^{d},\omega \in \mathbb{R}^{r}.
\end{align*}
\end{assump}
This assumption is often made in the existing literature on two-time-scale stochastic approximation \cite{hong2020two,doan2021finite} and is a sufficient condition to guarantee the fast convergence of the auxiliary variable iterate.
This assumption essentially states that $G$ behaves similarly to the gradient of a strongly convex function in expectation, though it may not even be a gradient mapping.
It can be verified that Assumption \ref{assump:stronglymonotone_G} (or an equivalent contractive condition) holds in the actor-critic methods discussed in Section \ref{sec:applications}. 

In addition, we assume that $\omega^{\star}(\cdot)$ is Lipschitz continuous with respect to $\theta$.    
\begin{assump}\label{assump:Lipschitz_omega}
There exists a constant $L,B>0$ such that
\begin{align*}
    \left\|\omega^{\star}(\theta)-\omega^{\star}(\theta')\right\| \leq L\|\theta-\theta'\|, \quad \left\|\omega^{\star}(\theta)\right\|\leq B, \quad \forall \theta,\theta' \in \mathbb{R}^{d}.
\end{align*}
\end{assump}

Given two probability distributions $\mu_1$ and $\mu_2$ over the space $\Xcal$, their total variation (TV) distance is defined as
\begin{equation}
    \label{eq:TV_def}
    d_{\text{TV}}(\mu_1,\mu_2)=\frac{1}{2} \sup _{\nu: \Xcal \rightarrow[-1,1]}\left|\int \nu d \mu_1-\int \nu d \mu_2\right|.
\end{equation}
The definition of the mixing time of a Markov chain $\{X_{k}\}$ is given as follows.
\begin{definition}
    \label{def:mixing_time}
    Consider the Markov chain $\{X_k^{\theta}\}$ generated according to $X_k^{\theta}\sim\Pcal(\cdot\mid X_{k-1}^{\theta},\theta)$, and let $\mu_{\theta}$ be its stationary distribution.  For any $\alpha>0$, the mixing time of the chain $\{X_k^{\theta}\}$ corresponding to $\alpha$ is defined as
    \[
        \tau_{\theta}(\alpha) = \min\{k\in\mathbb{N}:\sup_{X\in\Xcal}d_{\text{TV}}(P(X_k^{\theta}=\cdot\mid X_0^{\theta}=X),\mu_{\theta}(\cdot))\leq \alpha\}.
    \]
\end{definition}
The mixing time $\tau_{\theta}(\alpha)$ essentially measures the time needed for the Markov chain $\{X_{k}^{\theta}\}$ to approach its stationary distribution \cite{LevinPeresWilmer2006}. We next consider the following important assumption that guarantees that the Markov chain induced by any static $\theta$ ``mixes'' geometrically.
\begin{assump}\label{assump:markov-chain}
    Given any $\theta$, the Markov chain $\{X_{k}\}$ generated by $\Pcal(\cdot\mid\cdot,\theta)$ has a unique stationary distribution $\mu_{\theta}$ and is uniformly geometrically ergodic. In other words, there exist constants $m>0$ and $\rho\in (0,1)$ independent of $\theta$ such that
    \begin{equation*}
        \sup_{X\in\Xcal}d_{\text{TV}}(P(X_k=\cdot\,|X_0=X,\,\theta),\mu_{\theta}(\cdot))\leq m\rho^k \text{ for all } \theta\in\mathbb{R}^d\text{~and }k\geq 0. 
    \end{equation*}
\end{assump}

Denoting $\tau(\alpha)\triangleq\sup_{\theta\in\mathbb{R}^d}\tau_{\theta}(\alpha)$, this assumption implies that there exists a positive constant $C$ depending only on $\rho$ and $m$ such that 
\begin{equation}
    \tau(\alpha) \leq C\log\left(1/\alpha\right).\label{assump:mixing:tau}
\end{equation} 
Assumption \ref{assump:markov-chain} is again standard in the analysis of algorithms under time-varying Markovian samples \cite{zou2019finite,wu2020finite}. 

We also consider the following assumption on the ensemble of transition kernels.
\begin{assump}
    \label{assump:tv_bound}
    Given two distributions $d,\hat{d}$ over $\Xcal$ and parameters $\theta,\hat{\theta}\in\mathbb{R}^d$, we draw the samples according to
    $X\sim d, X'\sim\Pcal(\cdot\mid X,\theta)$ and $\hat{X}\sim\hat{d},\hat{X}'\sim \Pcal(\cdot\mid\hat{X},\hat{\theta})$.
    We assume that there exists a constant $L>0$ such that
    \begin{align}
        &d_{\text{TV}}(P(X'=\cdot), P(\hat{X}'=\cdot)) \leq d_{\text{TV}}(d, \hat{d}) + L\|\theta-\hat{\theta}\|.
        \label{assump:tv_bound:eq1}
    \end{align}
    In addition, we assume that the stationary distribution is Lipschitz in $\theta$
    \begin{align}
        d_{\text{TV}}(\mu_{\theta}, \mu_{\hat{\theta}})\leq L\|\theta-\hat{\theta}\|.
        \label{assump:tv_bound:eq2}
    \end{align}
\end{assump}

This assumption amounts to a regularity condition on the transition probability matrix $\Pcal(\cdot\mid\cdot,\theta)$ as a function of $\theta$, and has been shown to hold in the reinforcement learning setting (see, for example, \cite[Lemma A1]{wu2020finite}). Without any loss of generality, we use the same constant $L$ in Assumptions \ref{assump:HG_smooth}--\ref{assump:tv_bound} and assume $B\geq 1$. We define\vspace{-2pt}
\begin{equation}
D = \max\{L+\max_{X\in\Xcal}\|G(0,0,X)\|,\|\omega^*(0)\|,B\},  
\label{notation:D}
\end{equation}
which is a finite constant since $\Xcal$ is compact. A simple consequence of Assumption \ref{assump:HG_smooth} is that for all $\theta\in\mathbb{R}^d,\omega\in\mathbb{R}^r,X\in\Xcal$
\begin{gather}
    \|G(\theta,\omega,X)\|\leq D(\|\theta\|+\|\omega\|+1),\quad\text{and}\quad
    \|\omega^{\star}(\theta)\|\leq D(\|\theta\|+1).\label{eq:G_affinelybounded_D}
\end{gather}
Finally, we assume the optimal solution set $\{\theta^{\star}\hspace{-2pt}:f(\theta^{\star})\leq f(\theta),\,\forall \theta\in\mathbb{R}^{d}\}$ is non-empty.

\section{Finite-Time Complexity of Two-Time-Scale SGD}\label{sec:main_theorems}
This section presents the main results of this paper, which are the finite-time convergence of Algorithm \ref{Alg:two-time-scale-SGD} under three structural properties of the objective function, namely, strong convexity,
non-convexity with the P\L~condition, and general non-convexity.
Our results are derived under the assumptions introduced in Section \ref{sec:assumptions}, which we assume always hold in the rest of this paper.

The convergence of Algorithm \ref{Alg:two-time-scale-SGD} relies on $\alpha_k,\beta_k\rightarrow 0$ with reasonable rates. As mentioned in Section \ref{sec:two-time-scale-SGD}, $\alpha_{k}$ needs to be much smaller than $\beta_{k}$ to approximate the nested-loop algorithm where multiples auxiliary variable updates are performed for each decision variable update. Therefore, we consider the following choices of step sizes\looseness=-1
\begin{align}
    \alpha_k = \frac{\alpha_0}{(k+1)^{a}}, \quad \beta_k = \frac{\beta_0}{(k+1)^{b}},\quad\forall k\geq 0,\label{stepsizes}
\end{align}
where $a,b,\alpha_0,\beta_0$ are some constants satisfying $0 < b\leq a\leq 1$ and $0< \alpha_0 \leq  \beta_0$. 
Given $\alpha_{k}$, recall from Definition \ref{def:mixing_time} that  $\tau(\alpha_{k})$ is the mixing time associated with $\alpha_{k}$. In the sequel, for convenience we denote $\tau_{k}\triangleq\tau(\alpha_k)$.
Since $\tau_{k} \leq C\log((k+1)^{a}/\alpha_0)$ (from \eqref{assump:mixing:tau}), we have
$\lim_{k\rightarrow\infty}\alpha_{k}\tau_{k}^2 = \lim_{k\rightarrow\infty}\beta_{k}\tau_{k}^2 = 0$.
This implies that there exists a positive integer $\Kcal$ such that
\begin{align}
    \beta_{k-\tau_k}\tau_k^2\leq\min\Big\{1,\frac{1}{6LB},\frac{\lambda}{22C_1+32D^2},\frac{\lambda}{2C_2},\frac{\lambda^3}{32L^2 C_2}\Big\},\quad\forall k\geq\Kcal,\label{notation:Kcal}
\end{align}
where $C_1$ and $C_2$ are positive constants defined as 
\begin{align}
C_1 &= 18D^2+20LDB,\,\,C_2=(4D^2+1)(4C_1+32D^2)+\frac{2L^2 B^2}{\lambda}+2L^2 B^2.
\label{notation:constant_C}
\end{align}
In addition, there exists a constant $c_{\tau}\in(0,1)$ such that for any $k>\tau_k$ we have
\begin{align}
    \tau_k\leq(1-c_{\tau})k+(1-c_{\tau}), \quad\text{and}\quad c_{\tau}(k+1)\leq k-\tau_k+1\leq k+1.
\end{align}

We carefully come up with the constants and conditions in \eqref{notation:Kcal} and \eqref{notation:constant_C} to prevent an excessively large step size from destroying the stability of the updates. We note that $\Kcal$ is a constant that only depends on the quantities involved in the step sizes in \eqref{stepsizes}.

\subsection{Strong Convexity} We consider the following assumption on function $f$. 

\begin{assump}\label{assump:sc}
The function $f$ is strongly convex with constant $\lambda$\footnote{Without any loss of generality, we slightly overload $\lambda$, the strong monotonicity constant of the operator $G$ in Assumption~\ref{assump:stronglymonotone_G},
to denote the strong convexity constant here.}
\begin{align}
    f(y) \geq f(x)+\langle\nabla f(x),y-x\rangle+\frac{\lambda}{2}\|y-x\|^{2}, \quad \forall x,y\in\mathbb{R}^d.\label{assump:sc:ineq}
\end{align}
\end{assump}
\begin{thm}[Strongly Convex]
Suppose that Assumption \ref{assump:sc} holds.
Let the step size sequences $\{\alpha_k\}$ and $\{\beta_k\}$ satisfy \eqref{stepsizes} with
\begin{align*}
a=1,\quad b=2/3,\quad \alpha_0\geq \frac{4}{\lambda},\quad\text{and}\quad \frac{\alpha_0}{\beta_0}\leq\frac{1}{2}\cdot
\end{align*}
Then for all $k\geq \Kcal$ where $\Kcal$ is defined in \eqref{notation:Kcal},
we have
\begin{align*}
    \mathbb{E}\left[\|\theta_k-\theta^{\star}\|^2\right]
    &\leq \frac{\Kcal+1}{k+1}\left(\mathbb{E}\left[\|\theta_{\Kcal}-\theta^{\star}\|^2\right]+\frac{4 L^{2} \alpha_0}{\lambda^2 \beta_0}\mathbb{E}[\|\omega_{\Kcal}-\omega^{\star}(\theta_{\Kcal})\|^2]\right)\notag\\
    &\hspace{-43pt}+\frac{C^2\log^2((k+1)/\alpha_0)}{3(k+1)^{2/3}}\left((6C_2+2B^2+\frac{8L^2 C_2}{\lambda^2}(2\|\theta^{\star}\|^2+1))\frac{\alpha_0\beta_0}{c_{\tau}}+\frac{16L^4 B^2 \alpha_0^3}{\lambda^3 \beta_0^2}\right).
\end{align*}
\label{thm:convergence_decision_sc}
\end{thm}
Our theorem states that when $f$ is strongly convex the iterates of Algorithm~\ref{Alg:two-time-scale-SGD} converge to the optimal solution with rate $\widetilde{\Ocal}(k^{-2/3})$.
Comparing with the deterministic gradient descent setting where the convergence rate is linear and the standard SGD setting where the convergence rate is $\Ocal(k^{-1})$, our result reflects the compromise in the convergence rate due to the gradient noise and inaccurate auxiliary variable.
Compared with the convergence rate of the two-time-scale SA algorithm for bi-level optimization \cite{hong2020two} under i.i.d.\ samples, our rate is the same up to a logarithmic factor which naturally arises from the bias caused by the time-varying Markovian samples.


\subsection{Non-Convexity under P\L~Condition}\label{sec:conv_pl}
We also study the convergence of Algorithm \ref{Alg:two-time-scale-SGD} under the following condition.
\begin{assump}\label{assump:PL_condition}
There exists a constant $\lambda>0$ such that
\[\frac{1}{2}\|\nabla f(x)\|^2 \geq \lambda \left(f(x)-f^{\star}\right), \quad \forall x\in\mathbb{R}^d.\]
\end{assump}
This is known as the P\L~condition and is introduced in \cite{polyak1987introduction,lojasiewicz1963topological}. The P\L~condition does not imply convexity, but guarantees the linear convergence of the objective function value when gradient descent is applied to solve a non-convex optimization problem \cite{karimi2016linear}, which resembles the convergence rate of gradient descent for strongly convex functions. Recently, this condition has been observed to hold in many important practical problems such as supervised learning with an over-parametrized neural network \cite{liu2020toward} and the linear quadratic regulator in optimal control \cite{fazel2018global,yang2019global}.

\begin{thm}[P\L~Condition]\label{thm:convergence_decision_pl}
Suppose the function $f$ satisfies Assumption~\ref{assump:PL_condition}.
In addition, we assume that the stochastic gradient is bounded, i.e. there exists a constant $B>0$ such that \footnote{Again, for the convenience of notation, we use the same constant $B$ as in Assumption~\ref{assump:Lipschitz_omega}.}
\begin{align}
\|H(\theta,\omega,X)\| \leq B,\quad \forall \theta\in\Rset^{d}, \omega\in\Rset^{r},X\in\Xcal.    \label{assump:H_bounded}
\end{align}
Let the step size sequences $\{\alpha_k\}$ and $\{\beta_k\}$ satisfy \eqref{stepsizes} with
\begin{align*}
    a=1,\quad b=2/3,\quad\alpha_0\geq \max\{1,\frac{2}{\lambda}\},\quad\text{and}\quad\frac{\alpha_0}{\beta_0}\leq\frac{1}{4}.
\end{align*}
Then for all $k\geq \Kcal$ where $\Kcal$ is defined in \eqref{notation:Kcal}, we have
\begin{align*}
    \mathbb{E}\left[f(\theta_{k})\hspace{-2pt}-\hspace{-2pt}f^{\star}\right]&\hspace{-3pt}\leq\hspace{-3pt} \frac{\Kcal+1}{k+1}\hspace{-2pt}\left(\hspace{-2pt}\mathbb{E}\left[f(\theta_{\Kcal})\hspace{-2pt}-\hspace{-2pt}f^{\star}\right]\hspace{-2pt}+\hspace{-2pt}\frac{2L^{2} \alpha_0}{\lambda \beta_0}\mathbb{E}[\|\omega_{\Kcal}-\omega^{\star}(\theta_{\Kcal})\|^2]\hspace{-2pt}\right)\hspace{-2pt}+\hspace{-2pt}\frac{2C^2 C_3\log^4(k\hspace{-2pt}+\hspace{-2pt}1)}{3(k+1)^{2/3}},
\end{align*}
where $C_3=\frac{150L^2B^3\alpha_0^2}{c_{\tau}}+\frac{48L^2 C_2}{\lambda} (\|\theta_0\|^2+B^2 \alpha_0^2+1)\frac{\alpha_0\beta_0}{c_{\tau}}+\frac{48L^4 B^2 \alpha_0^3}{\lambda^2 \beta_0^2}$.
\end{thm}

Under the P\L~condition, we show that $f(\theta_k)$ converges to the optimal function value $f^{\star}$ with rate $\widetilde{\Ocal}(k^{-2/3})$. This is the same rate as if $f$ is strongly convex. However, in this case the convergence is measured in the function value, whereas under strong convexity the iterates $\theta_k$ converge to the unique global minimizer. The convergence rates of deterministic gradient descent and standard SGD under the P\L~condition also match those in the strongly convex case. To our best knowledge, functions exhibiting the P\L~condition have not been studied in the bi-level optimization framework.

\subsection{Non-Convexity} 
Finally, we study the case where the objective function $f$ is non-convex and smooth. In general, we cannot find an optimal solution and may only reach a stationary point. 
Analyzing the convergence without any convexity or P\L~condition is more challenging, and we need to make an additional assumption to ensure stability. 

\begin{assump}\label{assump:HG_smooth_nonconvex}
There exists a constant $L > 0$ such that
\begin{align*}
    \left\|G(\theta_1,\omega_1,X)-G(\theta_2,\omega_2,X)\right\| \leq L(\|\omega_1-\omega_2\|+1),\forall \theta_{1},\theta_{2}\in\mathbb{R}^d,\omega_{1},\omega_{2}\in\mathbb{R}^r,X\in\Xcal.
\end{align*}
\end{assump}
We note that this assumption holds in the actor-critic algorithm we discuss in Section \ref{sec:MDP_AC} where $G$ does not depend on $\theta$, as well as in problems where $G$ is bounded in $\theta$.

\begin{thm}[Non-convex]\label{thm:convergence_decision_nonconvex}
Let the step size $\{\alpha_k\}$ and $\{\beta_k\}$ satisfy \eqref{stepsizes} with
$a=3/5$, and $b=2/5$.
Under Assumptions \ref{assump:HG_smooth_nonconvex}, we have for all $k\geq\Kcal$
\begin{align*}
    \min_{t\leq k}\mathbb{E}[\|\nabla f(\theta_t)\|^2] &\leq \frac{4}{\alpha_0(k+1)^{2/5}}\mathbb{E}[f(\theta_{\Kcal})-f^{\star}]+\frac{4L^2}{\beta_0\lambda (k+1)^{2/5}}\mathbb{E}[\|\omega_{\Kcal}-\omega^{\star}(\theta_{\Kcal})\|^2]\notag\\
    &\hspace{15pt}+\left(\frac{25L^2B^3\alpha_0^2}{2c_{\tau}}+\frac{2L^4 B^2 \alpha_0^3}{\lambda \beta_0^2}+\frac{L^2 \alpha_0\beta_0 C_2}{c_{\tau}\lambda}\right)\frac{8\tau_k^2\log(k+1)}{5\log(2)c_{\tau}(k+1)^{2/5}}.
\end{align*}
\end{thm}

Our theorem in the non-convex case shows the convergence of the two-time-scale SGD algorithm to a stationary point of the objective function (measured by the squared norm of the gradient) with rate $\widetilde{\Ocal}(k^{-2/5})$.
One may contrast this with the convergence rate of deterministic gradient descent $\Ocal(k^{-1})$ and standard SGD $\Ocal(k^{-1/2})$ to see the cost of the gradient noise and auxiliary variable inaccuracy. Compared with the bi-level optimization algorithm under i.i.d.\ samples \cite{hong2020two}, our rate is again the same up to a logarithmic factor due to the time-varying Markovian samples.

\section{Analysis Decomposition and Proof of Main Theorem}
In this section, we briefly explain the main technical challenge in analyzing Algorithm \ref{Alg:two-time-scale-SGD}, which is the coupling between $\theta$, $\omega$,
and the time-varying Markovian samples. 
Our approach to the challenge is to properly ``decouple'' the variable updates so that we can handle them individually. 
Specifically, we first show under time-varying Markovian samples the convergence of the decision variable up to an error in the auxiliary variable (Section~\ref{sec:decision_var_conv}) and the reduction of the auxiliary variable error which hinges on the decision variable convergence (Section~\ref{sec:auxiliary_var_conv}), which essentially form a coupled dynamical system. In Section \ref{sec:lyapunov_lemma}, we introduce an important lemma that performs Lyapunov analysis on a coupled dynamical system of two inequalities. This lemma is a unified tool to analyze our algorithm under different function structures and may be of independent interest in the study of the finite-time performance of multiple-time-scale dynamical systems apart from those considered in this paper. Finally, we prove the theorem under the P\L~condition in Section~\ref{sec:proof_mainresults}. Strongly convex and general non-convex functions can be treated with similar analytical techniques, and the full details of their analyses can be found in the companion paper \cite{zeng2021two}.

When $f$ observes the P\L~condition, we show $\|f(\theta_k)-f^{\star}\|^2\rightarrow 0$. We frequently employ a few quantities for which we introduce the following shorthand notations
\begin{align}
\begin{aligned}
    z_k&\triangleq\omega_k-\omega^{\star}(\theta_k),\\
    \widebar{\Delta H_{k}} &\triangleq H(\theta_k,\omega_k,X_k)-H(\theta_k,\omega^{\star}(\theta_k),X_k),\\
    \Delta H_k &\triangleq H(\theta_k,\omega^{\star}(\theta_k),X_k)-\mathbb{E}_{\hat{X}\sim\mu_{\theta_k}}[H(\theta_k,\omega^{\star}(\theta_k),\hat{X})],\\
    \Delta G_k &\triangleq G(\theta_k,\omega_k,X_k)-\mathbb{E}_{\hat{X}\sim\mu_{\theta_k}}[G(\theta_k,\omega_k,\hat{X})].
    \end{aligned}
    \label{eq:Delta_z_def}
\end{align}
We can think of $\widebar{\Delta H_k}$ as the bias in the stochastic gradient due to the inaccurate auxiliary variable and $\Delta H_{k}$ and $\Delta G_{k}$ as the errors that the Markovian samples cause to $H$ and $G$.

\subsection{Decision Variable Convergence}\label{sec:decision_var_conv}
In this section, we derive a recursive formula for the iteration-wise decision variable convergence measured in $\mathbb{E}[f(\theta_{k})-f^{\star}]$. As a first step, we have from the update rule \eqref{alg:update_theta} and the $L$-smoothness of $f$
\begin{align}
    &f(\theta_{k+1})\leq f(\theta_{k})+\langle\nabla f(\theta_k),\theta_{k+1}-\theta_k\rangle+\frac{L}{2}\|\theta_{k+1}-\theta_k\|^2\notag\\
    &= f(\theta_{k})-\alpha_k\langle\nabla f(\theta_k),H(\theta_k,\omega_k,X_k)\rangle+\frac{L \alpha_k^2}{2}\|H(\theta_k,\omega_k,X_k)\|^2\notag\\
    &= f(\theta_{k})-\alpha_k\langle\nabla f(\theta_k),\mathbb{E}_{\hat{X}\sim\mu_{\theta_k}}[H(\theta_k,\omega^{\star}(\theta_k),\hat{X})]\rangle\notag\\
    &\hspace{20pt}-\alpha_k\left\langle\nabla f(\theta_k), \Delta H_k\hspace{-2pt}+\hspace{-2pt}\widebar{\Delta H_k}\right\rangle+\frac{L \alpha_k^2}{2}\|H(\theta_k,\omega_k, X_k)\|^2\notag\\
    &= f(\theta_{k})\hspace{-2pt}-\hspace{-2pt}\alpha_k\|\nabla f(\theta_k)\|^2\hspace{-2pt}-\hspace{-2pt}\alpha_k\left\langle\nabla f(\theta_k), \Delta H_k\hspace{-2pt}+\hspace{-2pt}\widebar{\Delta H_k}\right\rangle+\frac{L \alpha_k^2}{2}\|H(\theta_k,\omega_k, X_k)\|^2,
    \label{prop:pl_nonconvex:eq1}
\end{align}
where the last equality follows from \eqref{operator:H}, i.e. $\nabla f(\theta_k)=\mathbb{E}_{\hat{X}\sim \mu_{\theta_k}}[H(\theta_k,\omega^{\star}(\theta_k),\hat{X})]$. A key challenge to overcome is the time-varying Markovian randomness. If the samples were i.i.d.\ and the auxiliary variables were always solved perfectly, we would have $\mathbb{E}[\Delta H_k]=\mathbb{E}[\widebar{\Delta H_k}]=0$, reducing the problem to the one studied in the standard SGD. In the following lemma, we carefully treat the Markovian noise by leveraging the uniform geometric mixing time of the time-varying Markov chain and the Lipschitz condition of the state transition kernel.
\begin{lem}\label{lem:Gamma_bound}
    For any $k\geq \tau_k$, we have
    \begin{align*}
    \mathbb{E}\left[-\left\langle\nabla f(\theta_k),\Delta H_k\right\rangle\right] &\leq 12L^2B^3\tau_k^2\alpha_{k-\tau_k}.
    \end{align*}
\end{lem}
\begin{proof}
Our Markov process is a time-varying one (they depend on the iterates $\theta$). Therefore, one cannot directly utilize Assumption \ref{assump:markov-chain} to analyze the bias of $G$ in Algorithm \ref{Alg:two-time-scale-SGD} since the mixing time is defined for a fixed Markov chain (see Definition \ref{def:mixing_time}). To handle this difficulty, we introduce an auxiliary Markov chain $\{\widetilde{X}_{k}\}$ generated under the decision variable $\theta_{k-\tau_k}$ starting from $X_{k-\tau_{k}}$ as follows
\begin{align}
{X}_{k-\tau_k} \stackrel{\theta_{k-\tau_k}}{\longrightarrow}  \widetilde{X}_{k-\tau_k+1} \stackrel{\theta_{k-\tau_k}}{\longrightarrow}  \cdots \widetilde{X}_{k-1} \stackrel{\theta_{k-\tau_k}}{\longrightarrow}  \widetilde{X}_{k}.
\label{eq:def_distribution_d}
\end{align}
For clarity, we recall original the time-varying Markov processes $\{X_k\}$ generated by Algorithm~\ref{Alg:two-time-scale-SGD}
\begin{align*}
X_{k-\tau_k} \stackrel{\theta_{k-\tau_k+1}}{\longrightarrow} X_{k-\tau_k+1} \stackrel{\theta_{k-\tau_k+2}}{\longrightarrow} \cdots \stackrel{\theta_{k-1}}{\longrightarrow} X_{k-1} \stackrel{\theta_{k}}{\longrightarrow} X_{k}.
\end{align*}
Using the shorthand notation $y_k=\nabla f(\theta_k)$, we define the following quantities
\begin{align*}
    &T_1 = \mathbb{E}[\langle y_k-y_{k-\tau_k},\mathbb{E}_{\hat{X}\sim\mu_{\theta_k}}[H(\theta_k,\omega^{\star}(\theta_k),\hat{X})]-H(\theta_k,\omega^{\star}(\theta_k),X_k)\rangle]\notag\\
    &T_2 = \mathbb{E}[\langle y_{k-\tau_k},H(\theta_{k-\tau_k},\omega^{\star}(\theta_{k-\tau_k}),X_k)-H(\theta_k,\omega^{\star}(\theta_k),X_k)\rangle]\notag\\
    &T_3 = \mathbb{E}[\langle y_{k-\tau_k},H(\theta_{k-\tau_k},\omega^{\star}(\theta_{k-\tau_k}),\widetilde{X}_k)-H(\theta_{k-\tau_k},\omega^{\star}(\theta_{k-\tau_k}),X_k)\rangle]\notag\\
    &T_4 = \mathbb{E}[\langle y_{k-\tau_k},\mathbb{E}_{\bar{X}\sim\mu_{\theta_{k-\tau_k}}}\hspace{-2pt}[H(\theta_{k-\tau_k},\omega^{\star}(\theta_{k-\tau_k}),\bar{X})]-H(\theta_{k-\tau_k},\omega^{\star}(\theta_{k-\tau_k}),\widetilde{X}_k)\rangle]\notag\\
    &T_5 = \mathbb{E}[\langle y_{k-\tau_k},\mathbb{E}_{\hat{X}\sim\mu_{\theta_k}}\hspace{-2pt}[H(\theta_{k-\tau_k},\hspace{-1pt}\omega^{\star}(\theta_{k-\tau_k}),\hspace{-1pt}\hat{X})]\hspace{-2pt}-\hspace{-2pt}\mathbb{E}_{\bar{X}\sim\mu_{\theta_{k-\tau_k}}}\hspace{-4pt}[H(\theta_{k-\tau_k},\hspace{-1pt}\omega^{\star}(\theta_{k-\tau_k}),\hspace{-1pt}\bar{X})]\rangle]\notag\\
    &T_6 = \mathbb{E}[\langle y_{k-\tau_k},\mathbb{E}_{\hat{X}\sim\mu_{\theta_k}}\hspace{-2pt}[H(\theta_k,\omega^{\star}(\theta_k),\hat{X})]-\mathbb{E}_{\hat{X}\sim\mu_{\theta_{k}}}[H(\theta_{k-\tau_k},\omega^{\star}(\theta_{k-\tau_k}),\hat{X}]\rangle].
\end{align*}
It is easy to see that
\begin{align}
    -\mathbb{E}\left[\left\langle\nabla f(\theta_k),\Delta H_k\right\rangle\right]
    &=T_1+T_2+T_3+T_4+T_5+T_6.
    \label{eq:bound_Gamma_breakdown}
\end{align}
We analyze the terms of \eqref{eq:bound_Gamma_breakdown} individually. First, we treat $T_1$ using the boundedness of $H$ and the Lipschitz continuity of $\nabla f$
\begin{align}
T_1&\leq \mathbb{E}\big[\left\|y_k-y_{k-\tau_k}\right\|\left\|H(\theta_k,\omega^{\star}(\theta_k),X_k)-\mathbb{E}_{\hat{X}\sim\mu_{\theta_k}}[H(\theta_k,\omega^{\star}(\theta_k),\hat{X})]\right\|\big]\notag\\
&\leq L\mathbb{E}\big[\left\|\theta_k-\theta_{k-\tau_k}\right\|\big]]\cdot 2B\leq2B^2 L\tau_k\alpha_{k-\tau_k},
\label{lem:Gamma_bound:T1}
\end{align}
where the last inequality follows from
\[\|\theta_k-\theta_{k-\tau_k}\|\leq\sum_{t=k-\tau_k}^{k}\|\alpha_t H(\theta_t,\omega_t,X_t)\|\leq B\tau_k\alpha_{k-\tau_k}.\]
Similarly, for $T_2$ we have
\begin{align}
    T_2&\leq \mathbb{E}\big[\|y_{k-\tau_k}\|\|H(\theta_k,\omega^{\star}(\theta_k),X_k)-H(\theta_{k-\tau_k},\omega^{\star}(\theta_{k-\tau_k}),X_k)\|\big]\notag\\
    &\leq B\mathbb{E}\big[\|H(\theta_k,\omega^{\star}(\theta_k),X_k)-H(\theta_{k-\tau_k},\omega^{\star}(\theta_{k-\tau_k}),X_k)\|\big]\notag\\
    &\leq BL\mathbb{E}\big[\|\theta_k-\theta_{k-\tau_k}\|+\|\omega^{\star}(\theta_k)-\omega^{\star}(\theta_{k-\tau_k})\|\big]\notag\\
    &\leq BL(L+1)\mathbb{E}\big[\|\theta_k-\theta_{k-\tau_k}\|\big]\leq B^2 L(L+1)\tau_k\alpha_{k-\tau_k}.
    \label{lem:Gamma_bound:T2}
\end{align}
To analyze $T_3$, we utilize the law of total expectation: given $\Fcal\subseteq\Fcal'$ and a random variable $X$ we have $\Eset[X\mid\Fcal] = \Eset[\Eset[X\mid\Fcal']\mid\Fcal]$.   

Let $\Fcal_k$ be $\Fcal_k=\{X_0,\ldots,X_k,\theta_0,\ldots,\theta_k,\omega_0,...,\omega_k\}$,
and for convenience we denote
\[p_{k}(x) = P(X_{k} = x\mid\Fcal_{k-1})\quad  \text{and}\quad  \tilde{p}_{k}(x) = P(\widetilde{X}_{k} = x\mid\Fcal_{k-1}).\]
Then, we have
\begin{align*}
    &\mathbb{E}\big[\left\langle y_{k-\tau_k},H(\theta_{k-\tau_k},\omega^{\star}(\theta_{k-\tau_k}),\widetilde{X}_k)-H(\theta_{k-\tau_k},\omega^{\star}(\theta_{k-\tau_k}),X_k)\right\rangle\mid \Fcal_{k-\tau_k}\big]\notag\\
    &\leq\left\| y_{k-\tau_k}\right\|\|\mathbb{E}[H(\theta_{k-\tau_k},\omega^{\star}(\theta_{k-\tau_k}),X_k)-H(\theta_{k-\tau_k},\omega^{\star}(\theta_{k-\tau_k}),\widetilde{X}_k)\mid \Fcal_{k-\tau_k}]\|\notag\\
    &=\left\|y_{k-\tau_k}\right\|\big\|\mathbb{E}\big[\mathbb{E}\big[H(\theta_{k-\tau_k},\omega^{\star}(\theta_{k-\tau_k}),X_k)-H(\theta_{k-\tau_k},\omega^{\star}(\theta_{k-\tau_k}),\widetilde{X}_k)\hspace{-2pt}\mid\hspace{-2pt} \Fcal_{k-1}\big]\mid \Fcal_{k-\tau_k}\big]\big\|\notag\\
    &\leq B\mathbb{E}[\int_{\Xcal}H(\theta_{k-\tau_k},\omega^{\star}(\theta_{k-\tau_k}),x)(p_k(x)-\widetilde{p}_k(x)) dx \mid \Fcal_{k-\tau_k}]\notag\\
    &\leq 2B^2 \mathbb{E}[d_{TV}(p_k(\cdot), \widetilde{p}_k(\cdot))\mid \Fcal_{k-\tau_k}]\notag\\
    &\leq 2B^2 \mathbb{E}\big[d_{TV}(p_{k-1}(\cdot),\widetilde{p}_{k-1}(\cdot))+L\|\theta_k-\theta_{k-\tau_k}\|\mid \Fcal_{k-\tau_k}\big],
\end{align*}
where the second inequality uses the definition of the TV distance in \eqref{eq:TV_def}, and the last inequality is a result of Assumption \ref{assump:tv_bound}. Recursively applying this inequality and taking the expectation, we get
\begin{align}
    T_3
    &\leq 2B^2 L\sum_{t=k-\tau_k+1}^{k-1}\mathbb{E}[\|\theta_t-\theta_{k-\tau_k}\|]\leq 2B^3 L\tau_k^2\alpha_{k-\tau_k}.
    \label{lem:Gamma_bound:T3}
\end{align}
Similarly, to bound $T_4$, we again use the definition of TV distance
\begin{align*}
    &\mathbb{E}[\langle y_{k-\tau_k},\mathbb{E}_{\bar{X}\sim\mu_{\theta_{k-\tau_k}}}[H(\theta_{k-\tau_k},\omega^{\star}(\theta_{k-\tau_k}),\bar{X})]-H(\theta_{k-\tau_k},\omega^{\star}(\theta_{k-\tau_k}),\widetilde{X}_k)\rangle\mid\Fcal_{k-\tau_k}]\notag\\
    &\leq\left\|y_{k-\tau_k}\right\|\big\|\mathbb{E}[H(\theta_{k-\tau_k},\omega^{\star}(\theta_{k-\tau_k}),\widetilde{X}_k)-\mathbb{E}_{\bar{X}\sim\mu_{\theta_{k-\tau_k}}}\hspace{-2pt}[H(\theta_{k-\tau_k},\omega^{\star}(\theta_{k-\tau_k}),\bar{X})]\hspace{-2pt}\mid\hspace{-2pt}\Fcal_{k-\tau_k}]\big\|\notag\\
    &\leq B\cdot 2B \mathbb{E}[d_{TV}(\widetilde{p}_k(\cdot),\mu_{\theta_{k-\tau_k}})\mid \Fcal_{k-\tau_k}].
\end{align*}
Taking the expectation and using the definition of the mixing time \eqref{def:mixing_time},
\begin{align}
    T_4&\leq 2B^2 \mathbb{E}[d_{TV}(P(\widetilde{X}_k=\cdot),\mu_{\theta_{k-\tau_k}})]
    \leq 2B^2\alpha_k.
    \label{lem:Gamma_bound:T4}
\end{align}
We next consider $T_5$
\begin{align*}
    &\mathbb{E}[\langle y_{k-\tau_k},\hspace{-1pt}\mathbb{E}_{\hat{X}\sim\mu_{\theta_k}}\hspace{-2pt}[H(\theta_{k-\tau_k},\hspace{-1pt}\omega^{\star}\hspace{-1pt}(\theta_{k-\tau_k}),\hspace{-1pt}\hat{X})]\hspace{-2pt}-\hspace{-2pt}\mathbb{E}_{\bar{X}\sim\mu_{\theta_{k-\tau_k}}}\hspace{-6pt}[H(\theta_{k-\tau_k},\hspace{-1pt}\omega^{\star}\hspace{-1pt}(\theta_{k-\tau_k}),\hspace{-1pt}\bar{X})]\rangle\hspace{-2pt}\mid\hspace{-2pt}\Fcal_{k-\tau_k}]\notag\\
    &\leq B\|\mathbb{E}[\mathbb{E}_{\bar{X}\sim\mu_{\theta_{k-\tau_k}}}\hspace{-4pt}[H(\theta_{k-\tau_k},\omega^{\star}(\theta_{k-\tau_k}),\hspace{-1pt}\bar{X})]\hspace{-2pt}-\hspace{-2pt}\mathbb{E}_{\hat{X}\sim\mu_{\theta_k}}\hspace{-2pt}[H(\theta_{k-\tau_k},\omega^{\star}(\theta_{k-\tau_k}),\hspace{-1pt}\hat{X})]\hspace{-2pt}\mid\hspace{-2pt}\Fcal_{k-\tau_k}]\|\notag\\
    &\leq 2B^2\mathbb{E}[d_{TV}(\mu_{\theta_{k-\tau_k}},\mu_{\theta_k})\mid \Fcal_{k-\tau_k}],
\end{align*}
where the last inequality again comes from the definition of the TV distance in \eqref{eq:TV_def}. By \eqref{assump:tv_bound:eq2} in Assumption \ref{assump:tv_bound}, we have
\begin{align}
    T_5&\leq2B^2\mathbb{E}[d_{TV}(\mu_{\theta_{k-\tau_k}},\mu_{\theta_k})]\leq 2B^2 L\mathbb{E}[\|\theta_k-\theta_{k-\tau_k}\|]\leq 2B^3 L\tau_k\alpha_k.
    \label{lem:Gamma_bound:T5}
\end{align}
Finally, we bound $T_6$ using the boundedness of $\nabla f$ and the Lipschitz continuity of $H$
\begin{align*}
    &T_6\leq \mathbb{E}\big[\| y_{k-\tau_k}\|\|\mathbb{E}_{\hat{X}\sim\mu_{\theta_{k}}}\big[H(\theta_{k-\tau_k},\omega^{\star}(\theta_{k-\tau_k}),\hat{X})\big]-\mathbb{E}_{\hat{X}\sim\mu_{\theta_k}}\big[H(\theta_k,\omega^{\star}(\theta_k),\hat{X})\big]\big]\|\notag\\
    &\leq B L\mathbb{E}[\|\theta_{k-\tau_k}\hspace{-2pt}-\hspace{-2pt}\theta_k\|\hspace{-2pt}+\hspace{-2pt}\|\omega^{\star}(\theta_{k-\tau_k})\hspace{-2pt}-\hspace{-2pt}\omega^{\star}(\theta_k)\|]\leq 2L^2B \mathbb{E}[\|\theta_{k-\tau_k}\hspace{-2pt}-\hspace{-2pt}\theta_k\|]\leq 2L^2B^2\tau_k\alpha_{k-\tau_k}.
\end{align*}
The claimed result follows from plugging the bounds on $T_1$-$T_6$ into \eqref{eq:bound_Gamma_breakdown}.
\end{proof}

We can use the Lipschitz continuity of $H$ to study the error caused by $\widebar{\Delta H_k}$ and show that it can be bounded by the sum of $\|\theta_k-\theta^{\star}\|^2$ and $\|z_k\|^2$. This bound on $\widebar{\Delta H_k}$ together with the bound on $\Delta H_k$ established in Lemma~\ref{lem:Gamma_bound} leads to the following proposition, which states that $\mathbb{E}[f(\theta_{k})-f^{\star}]$ is sufficiently reduced in every iteration if the auxiliary variable error $z_k$ is controlled.
\begin{prop}\label{prop:pl_nonconvex}
Under Assumptions \ref{assump:HG_smooth}-\ref{assump:tv_bound}, we have for all $k\geq\Kcal$
\begin{align*}
    \mathbb{E}[f(\theta_{k+1})-f^{\star}]\leq (1\hspace{-2pt}-\hspace{-2pt}\lambda\alpha_k)\mathbb{E}[f(\theta_{k})-f^{\star}]+\frac{L^2\alpha_k}{2}\mathbb{E}\left[\|z_k\|^2\right]+\frac{25L^2B^3}{2}\tau_k^2\alpha_k\alpha_{k-\tau_k}.
\end{align*}
\end{prop}
\begin{proof}

By the Lipschitz condition of the operator $H$,
\begin{align*}
    -\mathbb{E}[\langle\nabla f(\theta_k),\widebar{\Delta H_k}\rangle] &\leq\frac{1}{2}\mathbb{E}[\|\nabla f(\theta_k)\|^2+\|H(\theta_k, \omega_k, X_k)-H(\theta_k, \omega^{\star}(\theta_k), X_k)\|^2]\notag\\
    &\leq\frac{1}{2}\mathbb{E}[\|\nabla f(\theta_k)\|^2]+\frac{L^2}{2}\mathbb{E}\left[\|z_k\|^2\right].
\end{align*}
Using this inequality along with Lemma \ref{lem:Gamma_bound} in \eqref{prop:pl_nonconvex:eq1}, we have for all $k\geq \tau_k$ 
\begin{align*}
    \mathbb{E}[f(\theta_{k+1})]&\leq \mathbb{E}[f(\theta_{k})]-\alpha_k\left\langle\nabla f(\theta_k), \Delta H_k\right\rangle-\alpha_k\mathbb{E}[\|\nabla f(\theta_k)\|^2]\notag\\
    &\hspace{20pt}-\alpha_k\mathbb{E}\left[\left\langle\nabla f(\theta_k), \widebar{\Delta H_k}\right\rangle\right]+\frac{L B^2\alpha_k^2}{2}\notag\\
    &\leq \mathbb{E}[f(\theta_{k})]+12L^2B^3\tau_k^2\alpha_k\alpha_{k-\tau_k}-\alpha_k\mathbb{E}[\|\nabla f(\theta_k)\|^2]\notag\\
    &\hspace{20pt}+ \frac{\alpha_k}{2}\mathbb{E}[\|\nabla f(\theta_k)\|^2]+\frac{L^2\alpha_k}{2}\mathbb{E}\left[\|z_k\|^2\right]+\frac{L B^2 \alpha_k^2}{2}\notag\\
    &\leq \mathbb{E}[f(\theta_{k})]-\frac{\alpha_k}{2}\mathbb{E}[\|\nabla f(\theta_k)\|^2]+\frac{L^2\alpha_k}{2}\mathbb{E}\left[\|z_k\|^2\right]+\frac{25L^2B^3}{2}\tau_k^2\alpha_k\alpha_{k-\tau_k}\notag\\
    &\leq \mathbb{E}[f(\theta_{k})]-\lambda\alpha_k\mathbb{E}\left[f(\theta_k)-f^{\star}\right]+\frac{L^2\alpha_k}{2}\mathbb{E}\left[\|z_k\|^2\right]+\frac{25L^2B^3}{2}\tau_k^2\alpha_k\alpha_{k-\tau_k},
\end{align*}
where the last inequality is due to the P\L~condition.
Subtracting $f^{\star}$ from both sides of the inequality leads to the claimed result.
\end{proof}

\subsection{Auxiliary Variable Convergence}\label{sec:auxiliary_var_conv}
In this section, we present and analyze the convergence of the auxiliary variable, summarized in the proposition below.\looseness=-1
\begin{prop}\label{prop:auxiliary}
Under Assumptions \ref{assump:HG_smooth}-\ref{assump:tv_bound}, we have for all $k\geq\Kcal$
\begin{align*}
    \mathbb{E}[\|z_{k+1}\|^2] &\leq (1\hspace{-2pt}-\hspace{-2pt}\frac{\lambda\beta_k}{2})\mathbb{E}[\|z_{k}\|^2]\hspace{-2pt}+\hspace{-2pt}C_2\tau_k^2\beta_{k-\tau_k}\beta_k\mathbb{E}\left[\|\theta_k\|^2\right]\hspace{-2pt}+\hspace{-2pt}\frac{2L^2 B^2 \alpha_k^2}{\lambda\beta_k}+C_2\tau_k^2\beta_{k-\tau_k}\beta_k.
\end{align*}
\end{prop}

Recall the auxiliary variable error $z_k$ defined in \eqref{eq:Delta_z_def}. Proposition~\ref{prop:auxiliary} establishes an iteration-wise reduction of this error in expectation in face of the drift of $\theta_k$. To prove this proposition, we introduce Lemma~\ref{lem:Lambda_bound} that bounds the error in the auxiliary variable caused by the Markovian samples. We skip the proof of this lemma due to its similarity to Lemma~\ref{lem:Gamma_bound} and refer interested readers to our companion paper \cite{zeng2021two} for the full proof details.
\begin{lem}\label{lem:Lambda_bound}
    Recall the definition of $C_1$ in \eqref{notation:constant_C}. For any $k\geq\tau_k$, we have
    \begin{align*}
    \mathbb{E}[\langle z_k, \Delta G_k\rangle]\leq C_1\tau_k^2\beta_{k-\tau_k}\mathbb{E}\left[\|z_{k-\tau_k}\|^2+\|\theta_k\|^2+\|\omega_k\|^2+1\right].
\end{align*}
\end{lem}

Analyzing Proposition~\ref{prop:auxiliary} requires properly controlling $\|\omega_k-\omega_{k-\tau_k}\|$, which we handle in the following lemma.

\begin{lem}\label{lem:omega_k-omega_k-tau}
For all $k\geq\tau_k$, we have 
\begin{align*}
\|\omega_k-\omega_{k-\tau_k}\|&\leq 3D\beta_{k-\tau_k}\tau_k\left(\|\omega_{k}\| +\|\theta_{k}\|+1\right).
\end{align*}
\end{lem}
\begin{proof}
As a result of \eqref{eq:G_affinelybounded_D}, for any $k\geq 0$
\begin{align}
    \|\omega_{k+1}\|\hspace{-2pt}-\hspace{-2pt}\|\omega_{k}\|&\leq \|\omega_{k+1}\hspace{-2pt}-\hspace{-2pt}\omega_{k}\|\hspace{-2pt}=\hspace{-2pt}\|\beta_k G(\theta_{k+1},\omega_k,X_k)\|\hspace{-2pt}\leq\hspace{-2pt} D\beta_k\left(\|\theta_{k+1}\|\hspace{-2pt}+\hspace{-2pt}\|\omega_k\|\hspace{-2pt}+\hspace{-2pt}1\right).
    \label{lem:omega_k-omega_k-tau:ineq0}
\end{align}
Define $h_k=\|\omega_{k}\|+\|\theta_{k+1}\|$. We have for all $k\geq 1$
\begin{align*}
    h_k&=\|\omega_{k-1}+\beta_{k-1}G(\theta_k,\omega_{k-1},X_{k-1})\|+\|\theta_k+\alpha_k H(\theta_k,\omega_k,X_k)\|\notag\\
    &\leq \|\omega_{k-1}\|+D\beta_{k-1}(\|\theta_{k}\|+\|\omega_{k-1}\|+1)+\|\theta_k\|+B\alpha_k\notag\\
    &\leq (1+D\beta_{k-1})h_{k-1}+(B+D)\beta_{k-1}
\end{align*}
where the second inequality follows from \eqref{eq:G_affinelybounded_D} and Assumption \ref{assump:HG_smooth}, and the last inequality is due to $\alpha_k\leq\beta_k$ and the fact that $\{\beta_k\}$ is a decaying sequence.

Since $1+x \leq e^{x}$ for all $x\geq 0$, we have for all $k\geq\tau_k$ and $k-\tau_k\leq t\leq k$
\begin{align*}
    h_t&\leq (1+D\beta_{t-1})h_{t-1}+(B+D)\beta_{t-1}\notag\\
    &\leq (1+D\beta_{k-\tau_k})^{\tau_k+t-k} h_{k-\tau_k}+(B+D)\beta_{k-\tau_k}\sum_{t'=k-\tau_k}^{t-1}(1+D\beta_{k-\tau_k})^{t-t'-1}\notag\\
    &\leq (1+D\beta_{k-\tau_k})^{\tau_k} h_{k-\tau_k}+(B+D)\beta_{k-\tau_k}\tau_k(1+D\beta_{k-\tau_k})^{\tau_k}\notag\\
    &\leq e^{D\beta_{k-\tau_k}\tau_k}h_{k-\tau_k}+(B+D)\beta_{k-\tau_k}\tau_k e^{D\beta_{k-\tau_k}\tau_k}\leq 2h_{k-\tau_k}+\frac{1}{3},
\end{align*}
where the last inequality follows from the step size $2(B+D)\beta_{k-\tau_k}\tau_k\leq\frac{1}{3}\leq\log(2)$.
Combining this inequality with \eqref{lem:omega_k-omega_k-tau:ineq0}, we have for all $k\geq\tau_k$
\begin{align*}
    &\|\omega_k-\omega_{k-\tau_k}\|\leq \sum_{t=k-\tau_k}^{k-1}\hspace{-3pt}\|\omega_{t+1}-\omega_t\| \leq D\hspace{-3pt}\sum_{t=k-\tau_k}^{k-1}\hspace{-3pt}\beta_{t}\left(h_t + 1\right) \leq D\beta_{k-\tau_k}\hspace{-3pt}\sum_{t=k-\tau_k}^{k-1}\hspace{-3pt}\big(2h_{k-\tau_k} + \frac{4}{3}\big) \notag\\
    &\leq\hspace{-2pt} 2D\beta_{k-\tau_k}\hspace{-2pt}\tau_k(\|\omega_{k-\tau_k}\| \hspace{-2pt}+\hspace{-2pt}\|\theta_{k-\tau_k+1}\|\hspace{-2.5pt}+\hspace{-2.5pt} \frac{2}{3}) \hspace{-2pt}\leq\hspace{-2pt} 2D\beta_{k-\tau_k}\hspace{-2pt}\tau_k(\|\omega_{k-\tau_k}\| \hspace{-2pt}+\hspace{-2pt}\|\theta_{k}\|\hspace{-2pt}+\hspace{-2pt}B\beta_{k-\tau_k}\hspace{-2pt}(\tau_k\hspace{-2pt}-\hspace{-2pt}1)\hspace{-2.5pt}+\hspace{-2.5pt} \frac{2}{3})\notag\\
    &\leq 2D\beta_{k-\tau_k}\tau_k\left(\|\omega_k-\omega_{k-\tau_k}\|+\|\omega_{k}\| +\|\theta_{k}\|+1\right).
\end{align*}
Re-arranging terms and again using the step size condition $2D\beta_{k-\tau_k}\tau_k\leq\frac{1}{3}$, we get
\begin{align*}
\|\omega_k-\omega_{k-\tau_k}\|&\leq 3D\beta_{k-\tau_k}\tau_k\left(\|\omega_{k}\| +\|\theta_{k}\|+1\right).
\end{align*}
\end{proof}

\begin{proof}[Proof (of Proposition~\ref{prop:auxiliary})]

Recall that $z_{k} = \omega_{k} - \omega^{\star}(\theta_{k})$. Using \eqref{alg:update_omega}, we have
\begin{align*}
    \|z_{k+1}\|^2
    &= \|\omega_k+\beta_k G(\theta_{k+1},\omega_k, X_k)-\omega^{\star}(\theta_{k+1})\|^2\notag\\
    &= \|\left(\omega_k-\omega^{\star}(\theta_k)\right)+\beta_k G(\theta_{k+1},\omega_k,X_k)+\left(\omega^{\star}(\theta_k)-\omega^{\star}(\theta_{k+1})\right)\|^2\notag\\
    &\leq \|z_k\|^2+2\beta_k\langle z_k,G(\theta_{k+1},\omega_k,X_k)\rangle+2\langle z_k,\omega^{\star}(\theta_k)-\omega^{\star}(\theta_{k+1})\rangle\notag\\
    &\qquad + 2\beta_k^2\left\|G(\theta_{k+1},\omega_k,X_k)\right\|^2+2\|\omega^{\star}(\theta^k)-\omega^{\star}(\theta^{k+1})\|^2.
\end{align*}

From the definition of $\Delta G_k$ in \eqref{eq:Delta_z_def},
\begin{align}
    \|z_{k+1}\|^2&\leq \|z_k\|^2+2\beta_k\langle z_k,G(\theta_{k+1},\omega_k,X_k)-G(\theta_{k},\omega_k,X_k)\rangle+2\beta_k\langle z_k,G(\theta_{k},\omega_k,X_k)\rangle\notag\\
    &\quad+2\langle z_k,\omega^{\star}(\theta_k)\hspace{-2pt}-\hspace{-2pt}\omega^{\star}(\theta_{k+1})\rangle\hspace{-2pt} + \hspace{-2pt}2\beta_k^2\left\|G(\theta_{k+1},\omega_k,X_k)\right\|^2\hspace{-2pt}+\hspace{-2pt}2\|\omega^{\star}(\theta_k)-\omega^{\star}(\theta_{k+1})\|^2\notag\\
    &\leq \|z_k\|^2+2\beta_k\langle z_k,G(\theta_{k+1},\omega_k,X_k)-G(\theta_{k},\omega_k,X_k)\rangle\rangle\notag\\
    &\quad + 2\beta_k\langle z_k,\mathbb{E}_{\hat{X}\sim\mu_{\theta_{k}}}[G(\theta_{k},\omega_k,\hat{X})]\rangle \hspace{-2pt}+\hspace{-2pt} 2\beta_k\langle z_k,\Delta G_{k}\rangle + 2\langle z_k,\omega^{\star}(\theta_k)\hspace{-2pt}-\hspace{-2pt}\omega^{\star}(\theta_{k+1})\rangle \notag\\
    &\quad+2\beta_k^2\left\|G(\theta_{k+1},\omega_k,X_k)\right\|^2+2L^2 B^2\alpha_k^2,
    \label{prop:auxiliary_eq1}
\end{align}
where the second inequality applies Assumption \ref{assump:HG_smooth} and \eqref{alg:update_theta}, i.e.,
\begin{equation}
 \|\omega^{\star}(\theta_k)-\omega^{\star}(\theta_{k+1})\|^2 \leq  L^2 \|\theta_k-\theta_{k+1}\|^2= L^2\|\alpha_k H(\theta_k,\omega_k,X_k)\|^2 \leq L^2B^2\alpha_k^2.
    \label{prop:auxiliary_eq1.5}
\end{equation}
We next analyze each term on the right-hand side of \eqref{prop:auxiliary_eq1}. First, using the relation $\langle 2v_1,v_2\rangle\leq c\|v_1\|^2+\frac{1}{c}\|v_2\|^2$ for any vectors $v_1,v_2$ and scalar $c>0$, we bound the second term of \eqref{prop:auxiliary_eq1}
\begin{align}
    &\langle z_k,\hspace{-1pt}G(\theta_{k+1},\omega_k,\hspace{-1pt}X_k)\hspace{-2pt}-\hspace{-2pt}G(\theta_{k},\omega_k,\hspace{-1pt}X_k)\rangle\hspace{-2pt}\leq\hspace{-2pt} \frac{\lambda}{4}\|z_k\|^2\hspace{-2pt}+\hspace{-2pt}\frac{1}{\lambda}\|G(\theta_{k+1},\omega_k,\hspace{-1.5pt}X_k)\hspace{-2pt}-\hspace{-2pt}G(\theta_{k},\omega_k,\hspace{-1.5pt}X_k)\|^2\notag\\
    &\leq \frac{\lambda}{4}\|z_k\|^2+\frac{L^2}{\lambda}\|\theta_{k+1}-\theta_{k}\|^2 \leq \frac{\lambda}{4}\|z_k\|^2+\frac{L^2 B^2\alpha_k^2}{\lambda},\label{prop:auxiliary_eq1a}
\end{align}
where the second inequality follows from the Lipschitz continuity of $G$ and the last inequality is due to \eqref{prop:auxiliary_eq1.5}. Similarly, we consider the fifth term of \eqref{prop:auxiliary_eq1}
\begin{align}
    &2\langle z_k,\omega^{\star}(\theta_k)-\omega^{\star}(\theta_{k+1})\rangle\leq \frac{\beta_k\lambda}{2}\|z_k\|^2+\frac{2}{\lambda\beta_k}\|\omega^{\star}(\theta_k)-\omega^{\star}(\theta_{k+1})\|^2\notag\\
    &\leq \frac{\beta_k\lambda}{2}\|z_k\|^2+\frac{2L^2}{\lambda\beta_k}\|\theta_k-\theta_{k+1}\|^2\leq \frac{\beta_k\lambda}{2}\|z_k\|^2+\frac{2L^2 B^2 \alpha_k^2}{\lambda\beta_k}.\label{prop:auxiliary_eq1b}
\end{align}
Next, using Assumption \ref{assump:stronglymonotone_G} and $z_{k} = \omega_{k} - \omega^{\star}(\theta_{k})$ we consider the third term of \eqref{prop:auxiliary_eq1}
\begin{align}
    2\beta_k\langle z_k,\mathbb{E}_{\hat{X}\sim\mu_{\theta_{k}}}[G(\theta_{k},\omega_k,\hat{X})]\rangle \leq -2\lambda\beta_k\|z_k\|^2.\label{prop:auxiliary_eq1c}
\end{align}
By \eqref{eq:G_affinelybounded_D}  and \eqref{alg:update_theta} we have
\begin{align}
\|G(\theta_{k+1},\omega_{k},X_{k})\|^2 \hspace{-2pt}\leq \hspace{-2pt}2D^2 \big(\|\theta_{k+1}\| \hspace{-2pt}+\hspace{-2pt} \|\omega_{k}\| + 1\big)^2\hspace{-2pt}\leq \hspace{-2pt}  2D^2 \big(\|\theta_{k}\| \hspace{-2pt}+\hspace{-2pt} B\alpha_{k} + \|\omega_{k}\| \hspace{-2pt}+\hspace{-2pt} 1\big)^2.\label{prop:auxiliary_eq1d}
\end{align}
Taking the expectation on both sides of \eqref{prop:auxiliary_eq1} and using \eqref{prop:auxiliary_eq1a}--\eqref{prop:auxiliary_eq1d} and Lemma \ref{lem:Lambda_bound}
\begin{align}
    \mathbb{E}[\|z_{k+1}\|^2]
    &\leq \mathbb{E}[\|z_k\|^2]+\frac{\beta_k\lambda}{2}\mathbb{E}[\|z_k\|^2]+\frac{2L^2 B^2 \beta_k\alpha_k^2}{\lambda}+\frac{\beta_k\lambda}{2}\mathbb{E}[\|z_k\|^2]+\frac{2L^2 B^2 \alpha_k^2}{\lambda\beta_k}\notag\\
    &\quad-2\lambda\beta_k\mathbb{E}[\|z_k\|^2]+C_1\tau_k^2\beta_{k-\tau_k}\beta_k\mathbb{E}\left[\|z_{k-\tau_k}\|^2+\|\theta_k\|^2+\|\omega_k\|^2+1\right]\notag\\
    &\quad + 2D^2\beta_{k}^2 \big(\|\theta_{k}\| + B\alpha_{k} + \|\omega_{k}\| + 1\big)^2 + 2L^2 B^2\alpha_k^2\notag\\
    &\leq (1-\lambda\beta_k)\mathbb{E}[\|z_k\|^2]+C_1\tau_k^2\beta_{k-\tau_k}\beta_k\mathbb{E}\left[\|z_{k-\tau_k}\|^2\right]\notag\\
    &\quad+(C_1+8D^2)\tau_k^2\beta_{k-\tau_k}\beta_k\mathbb{E}\left[\|\theta_k\|^2+\|\omega_k\|^2\right]\notag\\
    &\quad+\frac{2L^2 B^2 \alpha_k^2}{\lambda\beta_k}+(C_1+32D^2+\frac{2L^2 B^2}{\lambda}+2L^2 B^2)\tau_k^2\beta_{k-\tau_k}\beta_k,
    \label{prop:auxiliary_eq2}
\end{align}
where the last inequality uses $\alpha_k\leq\beta_k$ and $B\alpha_k\leq 1$.  Note that $\|z_{k-\tau_k}\|^2$ obeys
\begin{align*}
    \|z_{k-\tau_k}\|^2&=\|z_k-(\omega_k-\omega_{k-\tau_k})+(\omega^{\star}(\theta_k)-\omega^{\star}(\theta_{k-\tau_k}))\|^2\notag\\
    &\leq 3\left(\|z_k\|^2+\|\omega_k-\omega_{k-\tau_k}\|^2+\|\omega^{\star}(\theta_k)-\omega^{\star}(\theta_{k-\tau_k})\|^2\right)\notag\\
    &\leq 3\|z_k\|^2+\frac{9}{4}\left(\left\|\omega_{k}\right\|^{2}+\left\|\theta_{k}\right\|^{2}+1\right)+L^2\|\theta_{k}-\theta_{k-\tau_k}\|^2\notag\\
    &\leq 3\left(\|z_k\|^2+\|\omega_{k}\|^{2}+\|\theta_{k}\|^{2}\right)+L^2 B^2\tau_k^2\alpha_{k-\tau_k}^2+\frac{9}{4}\notag\\
    &\leq 3\left(\|z_k\|^2+\|\omega_{k}\|^{2}+\|\theta_{k}\|^{2}+1\right),
\end{align*}
where the second inequality is due to Lemma \ref{lem:omega_k-omega_k-tau} and the Lipschitz continuity of $\omega^{\star}$, and the last inequality follows from the step size condition $L B\tau_k\alpha_{k-\tau_k}\leq\frac{1}{6}$. Substituting the preceding relation into \eqref{prop:auxiliary_eq2}, we have for all $k\geq\tau_k$
\begin{align}
    \mathbb{E}[\|z_{k+1}\|^2] &\leq (1-\lambda\beta_k)\mathbb{E}[\|z_k\|^2]+3C_1\tau_k^2\beta_{k-\tau_k}\beta_k\mathbb{E}\left[\|z_{k}\|^2\right]\notag\\
    &\quad+(4C_1+8D^2)\tau_k^2\beta_{k-\tau_k}\beta_k\mathbb{E}\left[\|\theta_k\|^2+\|\omega_k\|^2\right]\notag\\
    &\quad+\frac{2L^2 B^2 \alpha_k^2}{\lambda\beta_k}+(4C_1+32D^2+\frac{2L^2 B^2}{\lambda}+2L^2 B^2)\tau_k^2\beta_{k-\tau_k}\beta_k\notag\\
    &\leq (1-\lambda\beta_k)\mathbb{E}[\|z_k\|^2]+(11C_1+16D^2)\tau_k^2\beta_{k-\tau_k}\beta_k\mathbb{E}\left[\|z_{k}\|^2\right]\notag\\
    &\quad+(4D^2+1)(4C_1+8D^2)\tau_k^2\beta_{k-\tau_k}\beta_k\mathbb{E}\left[\|\theta_k\|^2\right]+\frac{2L^2 B^2 \alpha_k^2}{\lambda\beta_k}\notag\\
    &\quad+((4D^2+1)(4C_1+32D^2)+\frac{2L^2 B^2}{\lambda}+2L^2 B^2)\tau_k^2\beta_{k-\tau_k}\beta_k,
    \label{prop:auxiliary_eq2.5}
\end{align}
where in the last inequality we use \eqref{eq:G_affinelybounded_D} to derive 
\begin{align*}
    \|\omega_k\|^2\hspace{-2pt}\leq \hspace{-2pt}2\|\omega_{k}\hspace{-2pt}-\hspace{-2pt}\omega^{\star}(\theta_{k})\|^2\hspace{-2pt}+\hspace{-2pt}2\|\omega^{\star}(\theta_{k})\|^2\hspace{-2pt}\leq\hspace{-2pt} 2\|z_k\|^2\hspace{-2pt}+\hspace{-2pt}2D^2(\|\theta_k\|\hspace{-2pt}+\hspace{-2pt}1)^2\hspace{-2pt}\leq\hspace{-2pt} 2\|z_k\|^2\hspace{-2pt}+\hspace{-2pt}4D^2(\|\theta_k\|^2\hspace{-2pt}+\hspace{-2pt}1).
\end{align*}
By the choice of the step size we have $(11C_1+16D^2)\tau_k^2\beta_{k-\tau_k}\leq\frac{\lambda}{2}$. Thus, using the constant $C_2$ defined in \eqref{notation:constant_C}, \eqref{prop:auxiliary_eq2.5} implies
\begin{align*}
    \mathbb{E}[\|z_{k+1}\|^2] \leq(1-\frac{\lambda\beta_k}{2})\mathbb{E}[\|z_k\|^2]+C_2\tau_k^2\beta_{k-\tau_k}\beta_k\mathbb{E}\left[\|\theta_k\|^2+1\right]+\frac{2L^2 B^2 \alpha_k^2}{\lambda\beta_k}.
\end{align*}
\end{proof}
Propositions~\ref{prop:pl_nonconvex} and \ref{prop:auxiliary} show that the convergence of the decision variable and the auxiliary variable forms a coupled dynamical system that evolves under two different rates. In the next section, we introduce a two-time-scale lemma that solves the system.


\subsection{Two-Time-Scale Lemma}
\label{sec:lyapunov_lemma}
Although we analyze the performance of our algorithm for different types of objective functions and with different convergence metrics, these analyses eventually reduce to the study of two coupled inequalities. The dynamics of these two inequalities happen on different time scales determined by the two step sizes used in our algorithm. 
In this section we present a general result, which we call the two-time-scale lemma, that characterizes the behavior of these coupled inequalities. {\color{blue} At a first glance, part of our result in the following lemma resembles the coupled supermartingale lemma in \cite{wang2017stochastic}, but as we have explained in Section~\ref{sec:related_works} -- Bi-Level and Composite Optimization, the aim of \cite{wang2017stochastic} is to solve stochastic composite optimization, which considers stochastic oracles completely different from ours.}


\begin{lem}\label{lem:lyapunov}
    Let $\{a_k,\, b_k,\,c_k,\,d_k,\,e_k,\,f_k\}$ be non-negative 
    sequences satisfying $\frac{a_{k+1}}{d_{k+1}}\leq\frac{a_{k}}{d_{k}} < 1$, for all $k\geq0$.
    Let $\{x_k\}, \{y_k\}$ be two non-negative sequences. We consider two settings on their dynamics. 
\begin{enumerate}[leftmargin = 4.5mm]
    \item  Suppose that $x_{k},y_{k}$ satisfy the following coupled inequalities
    \begin{align}
    x_{k+1}\leq(1-a_k)x_k+b_k y_k+c_k,\quad y_{k+1}\leq(1-d_k)y_k+e_k x_k+f_k.\label{lem:lyapunov:sc:def}
    \end{align}
    In addition, assume that there exists a constant $A\in\mathbb{R}$ such that  
    \begin{align}
    Aa_k-b_k-\frac{A a_k^2}{d_k}\geq 0\quad \text{and}\quad  \frac{A e_k}{d_{k}}\leq\frac{1}{2}, \quad \text{for all } k\geq0. \label{lem:lyapunov:sc:def3}
    \end{align}
    Then we have  for all $0 \leq \tau \leq k$
    \begin{align*}
    x_k&\leq (x_{\tau}+\frac{A a_{\tau} }{d_{\tau}}y_{\tau})\prod_{t=\tau}^{k-1}(1-\frac{a_t}{2})+\sum_{\ell=\tau}^{k-1}\Big(c_{\ell}+\frac{A a_{\ell} f_{\ell}}{d_{\ell}}\Big) \prod_{t=\ell+1}^{k-1}(1-\frac{a_t}{2}).
    \end{align*}
    
    \item  Suppose that $\{x_k,y_k\}$ satisfy the following coupled inequalities
    \begin{align}
    x_{k+1}\leq(1+a_k)x_k+b_k y_k+c_k,\quad y_{k+1}\leq(1-d_k)y_k+e_k x_k+f_k.\label{lem:lyapunov:def}
    \end{align}    
    $\{u_k\}$ is a non-negative sequence such that
    \begin{align}
        u_k\leq(1+a_k)x_k-x_{k+1}+b_k y_k+c_k,\label{lem:lyapunov:def3}
    \end{align}
    then we have for any $0 \leq \tau \leq k$
    \begin{align*}
    \sum_{t=\tau}^{k}u_t \leq \Big(1+\sum_{t=\tau}^{k}(a_t+\frac{b_t e_t}{d_t})e^{\sum_{t=\tau}^{k}(a_t+\frac{b_t e_t}{d_t})}\Big)\Big(x_{\tau}+\frac{b_{\tau}y_{\tau}}{d_{\tau}}+\sum_{t=\tau}^{k}(c_t+\frac{b_t f_t}{d_t})\Big).
\end{align*}
\end{enumerate}    
\end{lem}
\begin{proof}
Case 1) Consider $V_k=x_k+\frac{A a_k}{d_k}y_k$. From the second equation in \eqref{lem:lyapunov:sc:def},
\begin{align*}
    \frac{A a_{k+1}}{d_{k+1}}y_{k+1} &\leq \frac{A a_{k}}{d_{k}}y_{k+1} \leq \frac{A a_{k}}{d_{k}}\left((1-d_k)y_k+e_k x_k+f_k\right)\notag\\
    &= (1-a_k)\frac{A a_{k}}{d_{k}}y_k+(a_k-d_k)\frac{A a_{k}}{d_{k}}y_k+\frac{A a_{k}e_k x_k}{d_{k}}+\frac{A a_{k} f_k}{d_{k}}.
\end{align*}
Combining this with the first inequality of \eqref{lem:lyapunov:sc:def} yields
\begin{align*}
    V_{k+1}&=x_{k+1}+\frac{A a_{k+1}}{d_{k+1}}y_{k+1}\notag\\
    &\leq (1\hspace{-2pt}-\hspace{-2pt}a_k)x_k\hspace{-2pt}+\hspace{-2pt}b_k y_k\hspace{-2pt}+\hspace{-2pt}c_k\hspace{-2pt}+\hspace{-2pt}(1\hspace{-2pt}-\hspace{-2pt}a_k)\frac{A a_{k}}{d_{k}}y_k\hspace{-2pt}+\hspace{-2pt}(a_k\hspace{-2pt}-\hspace{-2pt}d_k)\frac{A a_{k}}{d_{k}}y_k\hspace{-2pt}+\hspace{-2pt}\frac{A a_{k}e_k x_k}{d_{k}}\hspace{-2pt}+\hspace{-2pt}\frac{A a_{k} f_k}{d_{k}}\notag\\
    &= (1-a_k)\left(x_k+\frac{A a_k y_k}{d_k}\right)+\left(\frac{A a_k^2}{d_k}-A a_k+b_k\right)y_k+c_k+\frac{A a_{k}e_k x_k}{d_{k}}+\frac{A a_{k} f_k}{d_{k}}\notag\\
    &\leq (1-a_k)V_k+\frac{a_k}{2}x_k+c_k+\frac{A a_{k} f_k}{d_{k}}\leq (1-\frac{a_k}{2})V_k+c_k+\frac{A a_{k} f_k}{d_{k}},
\end{align*}
where the second inequality follows from \eqref{lem:lyapunov:sc:def3}.
Applying this relation recursively,
\begin{align*}
    x_k&\leq V_k\leq V_{\tau}\prod_{t=\tau}^{k-1}(1-\frac{a_t}{2})+\sum_{\ell=\tau}^{k-1}\left(c_{\ell}+\frac{A a_{\ell} f_{\ell}}{d_{\ell}}\right) \prod_{t=\ell+1}^{k-1}(1-\frac{a_t}{2})\\
    &\leq (x_{\tau}+\frac{A a_{\tau} }{d_{\tau}}y_{\tau})\prod_{t=\tau}^{k-1}(1-\frac{a_t}{2})+\sum_{\ell=\tau}^{k-1}\left(c_{\ell}+\frac{A a_{\ell} f_{\ell}}{d_{\ell}}\right) \prod_{t=\ell+1}^{k-1}(1-\frac{a_t}{2}).
\end{align*}

Case 2) Re-arranging the second inequality of \eqref{lem:lyapunov:def} and multiplying by $\frac{b_k}{d_k}$, we get\looseness=-1
\begin{align*}
    b_k y_k\leq \frac{b_k}{d_k}y_k-\frac{b_k}{d_k}y_{k+1}+\frac{b_k e_k x_k}{d_k}+\frac{b_k f_k}{d_k}.
\end{align*}
Plugging this inequality into the first inequality of \eqref{lem:lyapunov:def} yields
\begin{align}
    x_{k+1}&\leq(1+a_k)x_k+c_k+\frac{b_k}{d_k}y_k-\frac{b_k}{d_k}y_{k+1}+\frac{b_k e_k x_k}{d_k}+\frac{b_k f_k}{d_k}\notag\\
    &\leq(1+g_k)x_k+\frac{b_k}{d_k}y_k-\frac{b_k}{d_k}y_{k+1}+c_k+\frac{b_k f_k}{d_k},
    \label{lem:lyapunov:eq3}
\end{align}
where we define $g_k=a_k+\frac{b_k e_k}{d_k}$.
Since $1+c\leq\exp(c)$ for any scalar $c>0$, we have
\begin{align}
    x_{k+1}
    &\leq\exp(g_k)x_k+\frac{b_k}{d_k}y_k-\frac{b_k}{d_k}y_{k+1}+c_k+\frac{b_k f_k}{d_k}\notag\\
    &\leq\exp(\sum_{t=\tau}^{k}g_t)x_{\tau}+\exp\big(\sum_{t=\tau}^{k}g_t\big)\sum_{t=\tau}^{k}\left(\frac{b_t}{d_t}(y_t-y_{t+1})+c_t+\frac{b_t f_t}{d_t}\right)\notag\\
    &\leq \exp(\sum_{t=\tau}^{k}g_t)\Big(x_{\tau}+\frac{b_{\tau}y_{\tau}}{d_{\tau}}+\sum_{t=\tau}^{k}(c_t+\frac{b_t f_t}{d_t})\Big),
    \label{lem:lyapunov:eq4}
\end{align}
where the second inequality applies the first inequality recursively. 
The inequalities \eqref{lem:lyapunov:eq3}, \eqref{lem:lyapunov:eq4}, and \eqref{lem:lyapunov:def3} together imply
\begin{align*}
    \sum_{t=\tau}^{k}u_t&\leq \sum_{t=\tau}^{k}(x_t-x_{t+1})+\big(\hspace{-2pt}\max_{\tau\leq t\leq k}x_t\big)\hspace{-2pt}\sum_{t=\tau}^{k}g_t+\hspace{-2pt}\sum_{t=\tau}^{k}\left(\frac{b_t}{d_t}(y_t-y_{t+1})+c_t+\frac{b_t f_t}{d_t}\right)\notag\\
    &\leq x_{\tau}\hspace{-1pt}+\hspace{-2pt}\sum_{t=\tau}^{k}g_t\exp\big(\sum_{t=\tau}^{k}g_t\big)\Big(x_{\tau}+\frac{b_{\tau}y_{\tau}}{d_{\tau}}+\sum_{t=\tau}^{k}(c_t+\frac{b_t f_t}{d_t})\Big)+\frac{b_{\tau}}{d_{\tau}}y_{\tau}+\sum_{t=\tau}^{k}(c_t+\frac{b_t f_t}{d_t})\notag\\
    &=\left(1+\sum_{t=\tau}^{k}(a_t+\frac{b_t e_t}{d_t})\exp(\sum_{t=\tau}^{k}(a_t+\frac{b_t e_t}{d_t}))\right)\left(x_{\tau}+\frac{b_{\tau}y_{\tau}}{d_{\tau}}+\sum_{t=\tau}^{k}(c_t+\frac{b_t f_t}{d_t})\right).
\end{align*}
\end{proof}

Lemma \ref{lem:lyapunov} studies the behavior of the two interacting sequences $\{x_k\}$ and $\{y_k\}$ that have generic structure. In our analysis, properly selected convergence metrics on $\theta_k$ and $\omega_k$ evolve as $x_k$ and $y_k$ above, respectively, according to \eqref{lem:lyapunov:sc:def} for strongly convex and P\L~functions and \eqref{lem:lyapunov:def} for non-convex functions, while the sequences $\{a_k,\, b_k,\,c_k,\,d_k,\,e_k,\,f_k\}$ are ratios and products of the step sizes $\{\alpha_k\}$ and $\{\beta_k\}$.\looseness=-1

\subsection{Proof of Main Results}\label{sec:proof_mainresults}

In this section, we present the proof of Theorem~\ref{thm:convergence_decision_pl} which considers functions observing the P\L~condition. 
The analyses of strongly convex and general non-convex functions use similar techniques: one needs to properly select a convergence metric according to the function structure, set up a step-wise decay of the convergence metric like Proposition~\ref{prop:pl_nonconvex} which forms a coupled dynamical system with Proposition~\ref{prop:auxiliary}, and apply the two-time-scale lemma introduced in Section~\ref{sec:lyapunov_lemma} to the coupled system.

From the analysis of the auxiliary variable in Proposition \ref{prop:auxiliary}, we have for all $k\geq\Kcal$
\begin{align*}
    \mathbb{E}[\|z_{k+1}\|^2] &\leq (1\hspace{-2pt}-\hspace{-2pt}\frac{\lambda\beta_k}{2})\mathbb{E}[\|z_{k}\|^2]\hspace{-2pt}+\hspace{-2pt}C_2\tau_k^2\beta_{k-\tau_k}\beta_k\mathbb{E}\left[\|\theta_k\|^2\right]\hspace{-2pt}+\hspace{-2pt}\frac{2L^2 B^2 \alpha_k^2}{\lambda\beta_k}+C_2\tau_k^2\beta_{k-\tau_k}\beta_k.
\end{align*}
Due to the boundedness of the operator $H$,
\begin{align*}
    \|\theta_k\|\leq \|\theta_0\|+\sum_{t=0}^{k-1}\|\theta_{t+1}-\theta_t\|\leq \|\theta_0\|+\sum_{t=0}^{k-1}\frac{B \alpha}{t+1}\leq\|\theta_0\|+\frac{B \alpha\log(k+1)}{\log(2)},
\end{align*}
where the last inequality follows from $\sum_{t=0}^{t'}\frac{1}{(t+1)^{u}}\leq \frac{\log(t'+2)}{\log(2)}$ for any $t'\geq 0$. This relation implies for any $k\geq 0$
\begin{align*}
    \|\theta_k\|^2\leq 2\|\theta_0\|^2+\frac{2B^2 \alpha^2 \log^2(k+1)}{\log^2(2)}\leq 24(\|\theta_0\|^2+B^2 \alpha^2)\log^2(k+1).
\end{align*}
Using this inequality in the bound on $\mathbb{E}[\|z_{k+1}\|^2]$, we have
\begin{align*}
    \mathbb{E}[\|z_{k+1}\|^2]&\leq (1\hspace{-2pt}-\hspace{-2pt}\frac{\lambda\beta_k}{2})\mathbb{E}[\|z_{k}\|^2]\hspace{-2pt}+\hspace{-2pt}C_2\tau_k^2\beta_{k-\tau_k}\beta_k\mathbb{E}\left[\|\theta_k\|^2\right]\hspace{-2pt}+\hspace{-2pt}\frac{2L^2 B^2 \alpha_k^2}{\lambda\beta_k}\hspace{-2pt}+\hspace{-2pt}C_2\tau_k^2\beta_{k-\tau_k}\beta_k\notag\\
    &\hspace{-20pt}\leq \hspace{-2pt}(1\hspace{-2pt}-\hspace{-2pt}\frac{\lambda\beta_k}{2})\mathbb{E}[\|z_{k}\|^2]\hspace{-2pt}+\hspace{-2pt}24C_2 (\|\theta_0\|^2\hspace{-2pt}+\hspace{-2pt}B^2 \alpha^2\hspace{-2pt}+\hspace{-2pt}1)\tau_k^2\beta_{k-\tau_k}\beta_k\log^2(k\hspace{-2pt}+\hspace{-2pt}1)\hspace{-2pt}+\hspace{-2pt}\frac{2L^2 B^2 \alpha_k^2}{\lambda\beta_k}.
\end{align*}

We can apply Lemma \ref{lem:lyapunov} case 1) to the result of Proposition~\ref{prop:pl_nonconvex} and the inequality above with $\tau=\Kcal$ and
\begin{align*}
    & x_k\hspace{-2pt}=\hspace{-2pt}\mathbb{E}[f(\theta_{k})\hspace{-2pt}-\hspace{-2pt}f^{\star}],\,\, y_k\hspace{-2pt}=\hspace{-2pt}\mathbb{E}\left[\|z_k\|^2\right],\,\, a_k\hspace{-2pt}=\hspace{-2pt}\lambda\alpha_k,\,\, b_k\hspace{-2pt}=\hspace{-2pt}\frac{L^2\alpha_k}{2},\,\, c_k\hspace{-2pt}=\frac{25L^2B^3}{2}\tau_k^2\alpha_k\alpha_{k-\tau_k},\\
    & d_k=\frac{\lambda \beta_k}{2},\,\,e_k=0,\,\, f_k=24C_2 (\|\theta_0\|^2\hspace{-2pt}+\hspace{-2pt}B^2 \alpha^2\hspace{-2pt}+\hspace{-2pt}1)\tau_k^2\beta_{k-\tau_k}\beta_k\log^2(k\hspace{-2pt}+\hspace{-2pt}1)+\frac{2L^2 B^2 \alpha_k^2}{\lambda\beta_k}.
\end{align*}
In this case, one can verify that we can choose $A=\frac{L^2}{\lambda}$ if the step size sequences satisfy $\frac{\alpha_k}{\beta_k}\leq\frac{1}{4}$.
As a result of Lemma \ref{lem:lyapunov} case 1), we have for all $k\geq\Kcal$\vspace{-8pt}
\begin{align*}
    \mathbb{E}\left[f(\theta_{k})\hspace{-2pt}-\hspace{-2pt}f^{\star}\right]\hspace{-2pt}&\leq\hspace{-2pt} \left(\mathbb{E}\big[f(\theta_{\Kcal})-f^{\star}\right]+\frac{2L^{2} \alpha_{\Kcal}}{\lambda \beta_{\Kcal}}\mathbb{E}[\|z_{\Kcal}\|^2]\big)\hspace{-2pt}\prod_{t=\Kcal}^{k-1}\hspace{-2pt}(1\hspace{-2pt}-\hspace{-2pt}\frac{\lambda\alpha_t}{2})\hspace{-2pt}+\hspace{-2pt}\sum_{\ell=\Kcal}^{k-1}\prod_{t=\ell+1}^{k-1}(1\hspace{-2pt}-\hspace{-2pt}\frac{\lambda \alpha_t}{2})\notag\\
    &\hspace{-42pt}\times\Big(\frac{25L^2B^3}{2}\tau_k^2\alpha_{\ell}\alpha_{\ell-\tau_k}\hspace{-4pt}+\hspace{-3pt}\frac{48C_2 L^2}{\lambda} (\|\theta_0\|^2\hspace{-3pt}+\hspace{-3pt}B^2\hspace{-1pt} \alpha^2\hspace{-3pt}+\hspace{-3pt}1)\tau_k^2\beta_{k-\tau_k}\alpha_k\hspace{-1pt}\log^2\hspace{-1pt}(k\hspace{-3pt}+\hspace{-3pt}1)\hspace{-2pt}+\hspace{-2pt}\frac{4L^4 B^2 \alpha_k^3}{\lambda^2\beta_k^2}\Big).
\end{align*}
Plugging in the step sizes to the second term, we have
\begin{align}
    \mathbb{E}\left[f(\theta_{k})-f^{\star}\right]&\leq \left(\mathbb{E}\left[f(\theta_{\Kcal})-f^{\star}\right]+\frac{2L^{2} \alpha_{\tau_k}}{\lambda \beta_{\tau_k}}\mathbb{E}[\|z_{\Kcal}\|^2]\right)\prod_{t=\Kcal}^{k-1}(1-\frac{\lambda \alpha_t}{2})\notag\\
    &\hspace{10pt}+\tau_k^2\sum_{\ell=\Kcal}^{k-1}\Big(\frac{25L^2B^3\alpha_0^2}{2c_{\tau}(\ell+1)^2}+\frac{4L^4 B^2 \alpha_0^3}{\lambda^2\beta_0^2 (\ell+1)^{5/3}}\notag\\
    &\hspace{20pt}+\hspace{-2pt}\frac{48C_2 L^2}{\lambda} (\|\theta_0\|^2\hspace{-2pt}+\hspace{-2pt}B^2 \alpha_0^2\hspace{-2pt}+\hspace{-2pt}1)\frac{\alpha_0\beta_0}{c_{\tau}(\ell+1)^{5/3}}\log^2(k\hspace{-2pt}+\hspace{-2pt}1)\Big)\hspace{-4pt} \prod_{t=\ell+1}^{k-1}\hspace{-4pt}(1\hspace{-2pt}-\hspace{-2pt}\frac{\lambda\alpha_t}{2})\notag\\
    &\leq \left(\mathbb{E}\left[f(\theta_{\Kcal})-f^{\star}\right]+\frac{2L^{2} \alpha_{\Kcal}}{\lambda\beta_{\Kcal}}\mathbb{E}[\|z_{\Kcal}\|^2]\right)\prod_{t=\Kcal}^{k-1}(1-\frac{\lambda \alpha_t}{2})\notag\\
    &\hspace{10pt}+\tau_k^2\log^2(k+1)\sum_{\ell=\Kcal}^{k-1}\frac{C_3}{(\ell+1)^{5/3}} \prod_{t=\ell+1}^{k-1}(1-\frac{\lambda \alpha_t}{2}),
    \label{thm:convergence_decision_pl:eq2.6}
\end{align}
where we use the fact that $\frac{1}{\log^2(k+1)}\leq 12$ for all $k>0$ and the definition of $C_3$.

Since $1+c\leq\exp(c)$ for any scalar $c$, we have
\begin{align}
    \prod_{t=\Kcal}^{k-1}(1-\frac{\lambda \alpha_t}{2})&\leq \prod_{t=\Kcal}^{k-1}\exp(-\frac{\lambda \alpha_t}{2})=\exp(-\sum_{t=\Kcal}^{k-1}\frac{\lambda \alpha_t}{2})\leq \exp(-\frac{\lambda \alpha_0}{2}\sum_{t=\Kcal}^{k-1}\frac{1}{t+1})\notag\\
    &\leq \exp(-\frac{\lambda\alpha_0}{2}\log(\frac{k+1}{\Kcal+1}))\leq(\frac{\Kcal+1}{k+1})^{\frac{\lambda \alpha_0}{2}}\leq \frac{\Kcal+1}{k+1},
    \label{thm:convergence_decision_pl:eq2.7}
\end{align}
where the last inequality results from $\alpha_0\geq\frac{2}{\lambda}$, and the third inequality follows from $\sum_{t=k_1}^{k_2}\frac{1}{t+1}\geq \log(\frac{k_2+2}{k_1+1})$.
Similarly, we have
\begin{align}
    \prod_{t=\ell+1}^{k-1}(1-\frac{\lambda \alpha_t}{2})\leq \frac{2\ell+1}{k+1}\leq \frac{2(\ell+1)}{k+1}.
    \label{thm:convergence_decision_pl:eq2.8}
\end{align}
Using \eqref{thm:convergence_decision_pl:eq2.7} and \eqref{thm:convergence_decision_pl:eq2.8} in \eqref{thm:convergence_decision_pl:eq2.6}, \begin{align*}
    \mathbb{E}\left[f(\theta_{k})-f^{\star}\right]&\leq \left(\hspace{-2pt}\mathbb{E}\left[f(\theta_{\Kcal})\hspace{-2pt}-\hspace{-2pt}f^{\star}\right]\hspace{-2pt}+\hspace{-2pt}\frac{2L^{2} \alpha_{\Kcal}}{\lambda\beta_{\Kcal}}\mathbb{E}[\|z_{\Kcal}\|^2]\hspace{-2pt}\right)\hspace{-2pt}\frac{\Kcal}{k\hspace{-2pt}+\hspace{-2pt}1}\hspace{-2pt}+\hspace{-2pt}\frac{\log^2(k\hspace{-2pt}+\hspace{-2pt}1)\tau_k^2}{k+1} \hspace{-3pt}\sum_{\ell=\Kcal}^{k-1}\hspace{-4pt}\frac{2C_3}{(\ell\hspace{-2pt}+\hspace{-2pt}1)^{2/3}}\notag\\
    &\leq\frac{\Kcal+1}{k+1}\left(\mathbb{E}\left[f(\theta_{\Kcal})-f^{\star}\right]+\frac{2L^{2} \alpha_0}{\lambda \beta_0}\mathbb{E}[\|z_{\Kcal}\|^2]\right)+\frac{2C_3\log^2(k+1)\tau_k^2}{3(k+1)^{2/3}},
\end{align*}
where the second inequality is a result of the relation $\sum_{t=0}^{t'}\frac{1}{(t+1)^{2/3}}\leq\frac{(t'+1)^{1/3}}{3}$ for any $t'\geq0$. The claimed result follows from this and \eqref{assump:mixing:tau}.


\section{Details on Actor-Critic Algorithm for LQR}\label{sec:lqr_details}

This section provides further details on the policy optimization problem in LQR and the online actor-critic algorithm.

We define the matrix $\Sigma_K$ which helps us express the gradient of $J(K)$.  For fixed $K$, we define $\Sigma_K$ implicitly as the matrix that obeys the Lyapunov equation\looseness=-1
\begin{align}
    \Sigma_{K}=\Psi_{\sigma}+(A-B K) \Sigma_{K}(A-B K)^{\top}.
    \label{eq:SigmaK_def}
\end{align}
If $K$ is a stable feedback gain matrix, meaning that $\|A-BK\|<1$, then there exists a unique symmetric positive-definite $\Sigma_{K}$ observing this condition.  We can interpret $\Sigma_K$ as the covariance of the stationary distribution of the states under the feedback gain matrix $K$. In other words, the linear dynamical system in \eqref{eq:obj_LQR} has the stationary distribution $x\sim \mu_{K} = N(0,\Sigma_{K})$.  We can write the gradient of $J(K)$ as
\begin{align*}
    \nabla J(K) = E_{K} \Sigma_{K}.
\end{align*}

$J(K)$ is non-convex with respect to $K$ but satisfies the P\L~condition \cite[Lemma 5]{yang2019global}, which we will exploit for accelerated convergence,
\begin{align}
    &J(K)-J(K^{\star}) \leq \frac{\|\Sigma_{K^{\star}}\|}{\sigma_{\min }(\Sigma_{K})^2 \sigma_{\min }(R)}\|\nabla J(K)\|^2,
    \label{eq:lqr_PLcondition}
\end{align}
where $\sigma_{\min}(\cdot)$ returns the smallest singular value. Since $R$ is positive definite, $\sigma_{\min}(R)$ is a strictly positive constant. $\sigma_{\min}(\Sigma_K)$ is bounded away from zero under a proper assumption on $K$ as will be discussed below in Lemma \ref{lem:bounded_sigma_SigmaK}.

Estimating the gradient $\nabla J(K)$ requires estimating the two matrices $\Sigma_K$ and $E_K$, both of which depend on the unknown $A$ and $B$. While $E_K$ can be decomposed and estimated component-wise in a tractable manner using sampled data as we show in Section~\ref{sec:applications_lqr}, estimating $\Sigma_K$ is more difficult.
In light of this, our actor-critic algorithm will be based on estimating the natural policy gradient $\tilde{\nabla} J(K)$.
Natural policy gradient descent can be interpreted as policy gradient descent on a changed system of coordinates \cite{kakade2001natural}, and its analysis here will require only some small modifications to our general framework.

Detailed in Section~\ref{sec:proof_mainresults}, our analysis of gradient descent using the  P\L~condition centers on a second-order expansion of the objective function of the form 
\begin{align*}
    J(K_k)\leq J(K_{k-1})+\langle\nabla J(K_{k-1}), K_k-K_{k-1}\rangle+e_k.
\end{align*}
With gradient descent, we have $K_k-K_{k-1} = -\alpha_k\nabla J(K_{k-1})$ and the inner product term conveniently becomes $-\alpha_k\|\nabla J(K_{k-1})\|^2$, allowing us to eventually quantify the progress using \eqref{eq:lqr_PLcondition}. For natural gradient descent, we have $K_k-K_{k-1} = -\alpha_k\widetilde{\nabla} J(K_{k-1})$, and we can bound the inner product term as
\begin{align*}
    -\alpha_k\langle\nabla J(K_{k-1}),\widetilde{\nabla} J(K_{k-1})\rangle\leq-\alpha_k\sigma_{\min}(\Sigma_{K_{k-1}}^{-1})\|\nabla J(K_{k-1})\|^2.
\end{align*}
Since $\Sigma_K$ has upper and lower bounded eigenvalues as we will show later, we can absorb this multiplicative factor into the step size. 

\begin{algorithm}[t]
\caption{Online Actor Critic Algorithm for LQR}
\label{Alg:AC_LQR}
\begin{algorithmic}[1]
\STATE{\textbf{Initialization:} The control gain matrix $K_0\in\mathbb{R}^{d_1\times d_2}$, the auxiliary variable $\hat{\Omega}_0\in\mathbb{R}^{(d_1+d_2)\times(d_1+d_2)}, \hat{J}_0\in\mathbb{R}$, step size sequences $\{\alpha_k\}$ for the $K$ update, $\{\beta_k\}$ for the auxiliary variable update}
\STATE{Sample an initial state $x_0\sim\Dcal$, takes an initial control $u_0\sim N(-K_0 x_0,\sigma^2 I)$, and observe the cost $c_0=x_0^{\top}Qx_0+u_0^{\top}Ru_0$ and the next state $x_1\sim N(Ax_0+Bu_0,\Psi)$}
\FOR{$k=0,1,2,...$}
\STATE{Take the control $u_{k+1}\sim N(-K_{k} x_{k+1},\sigma^2 I)$ and observe the cost $c_{k+1}=x_{k+1}^{\top}Qx_{k+1}+u_{k+1}^{\top}Ru_{k+1}$ and the next state $x_{k+2}\sim N(Ax_{k+1}+Bu_{k+1},\Psi)$}
\STATE{\vspace{-10pt}Gain matrix update at the actor:
        \vspace{-5pt}
        \begin{align}
        K_{k+1} = K_k-\alpha_k(\hat{\Omega}_k^{22}K_k-\hat{\Omega}_k^{21})\label{Alg:AC_LQR:actor}
        \end{align}}
\STATE{\vspace{-15pt}Critic variable update:
        \vspace{-5pt}
        \begin{align}
        \hat{c}_k &= x_k^{\top} Q x_k + u_k^{\top} R u_k\notag\\
        \hat{J}_{k+1} &= \hat{J}_k - \beta_k(\hat{J}_k-\hat{c}_k)\label{Alg:AC_LQR:critic1} \\
        \hat{\Omega}_{k+1} &=\hat{\Omega}_{k} - \beta_k\left[\begin{array}{c}
        x_k \\
        u_k
        \end{array}\right]\left[\begin{array}{c}
        x_k \\
        u_k
        \end{array}\right]^{\top}\left(\left[
        x_k^{\top}\,\,\, u_k^{\top}\right]\hat{\Omega}_k\left[\begin{array}{c}
        x_k \\
        u_k
        \end{array}\right]+\hat{J}_k-\hat{c}_k\right)\label{Alg:AC_LQR:critic2}
        \end{align}
        }
\ENDFOR
\end{algorithmic}
\end{algorithm}

Our online actor-critic method for solving the LQR is formally presented in Algorithm \ref{Alg:AC_LQR}, where the actor $K_k$ and critic $[\hat{J}_k,\text{svec}(\hat{\Omega}_k)^{\top}]^{\top}$ are estimates of $K^{\star}$ and $[J(K_k),\text{svec}(\Omega_{K_k})^{\top}]^{\top}$ and are equivalent to $\theta_k$ and $\omega_k$ of Algorithm \ref{Alg:two-time-scale-SGD}. 
The purpose of the critic updates \eqref{Alg:AC_LQR:critic1} and \eqref{Alg:AC_LQR:critic2} is to find the solution to \eqref{eq:lqr_Goperator} with stochastic approximation, while \eqref{Alg:AC_LQR:actor} updates the actor along the natural gradient direction estimated by the critic $\hat{\Omega}_k$. 
As the problem obeys the P\L~condition \eqref{eq:lqr_PLcondition}, one can expect Algorithm \ref{Alg:AC_LQR} to converge in objective function value at a rate of $\Ocal(k^{-2/3})$ if Assumptions \ref{assump:HG_smooth}-\ref{assump:tv_bound} of the general optimization framework are satisfied. 
\begin{assump}
    \label{assump:stable_K}
    There exists a constant $\rho<1$ such that $\|A-BK_k\|< \rho,\forall k$ a.s.
\end{assump}
We make the assumption above on the uniform stability of $K_k$. Due to the special structure of the LQR, Assumption \ref{assump:stable_K} implies the boundedness of the gradient $\nabla J(K)$ and the natural gradient $\widetilde{\nabla} J(K)$ for all $K$ and the Lipschitz continuity of the operators $H$ and $G$, which translate to Assumptions \ref{assump:HG_smooth}, \ref{assump:f_smooth}, and \ref{assump:Lipschitz_omega}. In addition, \cite[Lemma 3.2]{yang2019global} shows that under Assumption \ref{assump:stable_K}, the matrix on the left hand side of \eqref{eq:lqr_Goperator} is positive definite, which implies the strong monotonicity of the operator $G(\theta,\omega,X)$, thus meeting Assumption \ref{assump:stronglymonotone_G}.
%
This assumption also implies Assumptions \ref{assump:markov-chain}-\ref{assump:tv_bound} and the boundedness of the singular values of $\Sigma_{K_k}$.
\begin{lem}
    \label{lem:bounded_sigma_SigmaK}
    Under Assumption \ref{assump:stable_K}, there exist constants $0<\sigma_l,\sigma_u<\infty$ depending only on $\Psi_{\sigma}$ and $\rho$ such that the eigenvalues of $\Sigma_{K_k}$ all lie within $[\sigma_l,\sigma_u]$ for any $k\geq0$.
\end{lem}
\begin{proof}
Given a matrix $M\in\mathbb{R}^{p\times q}$, we define $\lambda_{\max}(M)=\max_{x\in\mathbb{R}^q:\|x\|=1}\|Mx\|$ and $\lambda_{\min}(M)=\min_{x\in\mathbb{R}^q:\|x\|=1}\|Mx\|$.
From \eqref{eq:SigmaK_def}, for any $K$ such that $\|A-BK\|<\rho$\vspace{-3pt}
\begin{align*}
    &\lambda_{\max}(\Sigma_K)=\lambda_{\max}(\Psi_{\sigma}+(A-B K) \Sigma_{K}(A-B K)^{\top})\\
    &\leq \lambda_{\max}(\Psi_{\sigma})+\lambda_{\max}(A-B K)^2 \lambda_{\max}(\Sigma_K)\leq \lambda_{\max}(\Psi_{\sigma})+\rho^2 \lambda_{\max}(\Sigma_K),
\end{align*}
which yields $\lambda_{\max}(\Sigma_K) \leq \frac{\lambda_{\max}(\Psi_{\sigma})}{1-\rho^2}$. Assumption \ref{assump:stable_K} states that $\|A-BK_k\|<\rho$ for all $k$, implying that the largest singular value of the always positive definite matrix $\Sigma_{K_k}$, which is equal to $\lambda_{\max}(\Sigma_{K_k})$, is upper bounded.
On the other hand,\vspace{-3pt}
\begin{align*}
    \lambda_{\min}(\Sigma_K)&=\lambda_{\min}(\Psi_{\sigma}+(A-B K) \Sigma_{K}(A-B K)^{\top})\\
    &\geq \lambda_{\min}(\Psi_{\sigma})+\lambda_{\min}((A-B K) \Sigma_K (A-B K)^{\top}).
\end{align*}
Since $\Sigma_K$ is positive definite, $\lambda_{\min}((A-B K) \Sigma_K (A-B K)^{\top})\geq 0$, which implies $\lambda_{\min}(\Sigma_K)\geq \lambda_{\min}(\Psi_{\sigma})$.
\end{proof}

Lemma \ref{lem:bounded_sigma_SigmaK} ensures that the P\L~constant in \eqref{eq:lqr_PLcondition} is indeed upper bounded.

As we have verified the validity of all required assumptions, we can now state the following corollary on the convergence of Algorithm \ref{Alg:AC_LQR}.



\begin{cor}\label{cor:LQR}
Consider the iterates from $\Kcal$ iterations of Algorithm \ref{Alg:AC_LQR}.
Under Assumption \ref{assump:stable_K} and proper choices of $\alpha_0$ and $\beta_0$, we have for all $k\geq\Kcal$\vspace{-3pt}
\begin{align*}
\mathbb{E}\left[J(K_{k})-J(K^{\star})\right] &\leq \Ocal\Big(\frac{ \log ^{4}(k+1)}{(k+1)^{2 / 3}}\Big)\notag\\
&\hspace{-20pt}+\Ocal\Big(\frac{\Kcal+1}{k+1}\mathbb{E}\big[J(K_{\Kcal})-J(K^{\star})+\big\|\hat{\Omega}_{\Kcal}-\Omega_{K_\Kcal}\big\|^{2}+\big\|\hat{J}_{\Kcal}-J(K_\Kcal)\big\|^{2}\big]\Big).
\end{align*}
\end{cor}

This result guarantees that after a ``burn-in'' period that only depends on the mixing time of the Markov chain, the iterates $J(K_k)$ converges to the globally optimal cost of the re-formulated LQR problem \eqref{eq:obj_K_LQR} with a rate of $\widetilde{\Ocal}(k^{-2/3})$. 

\begin{figure}[ht]
  \includegraphics[width=\linewidth]{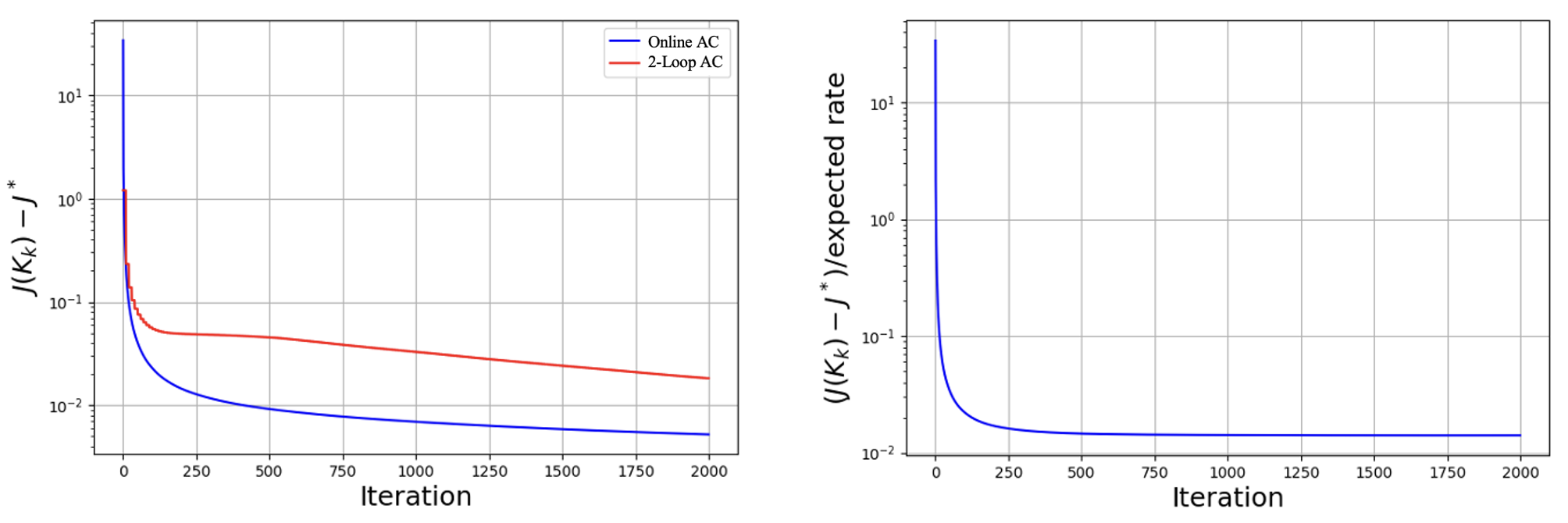}
  \vspace{-.5cm}
  \caption{Convergence rate of $J(K_k)-J^{\star}$ under Algorithm \ref{Alg:AC_LQR}}
  \label{fig:lqr_rate}
\end{figure}

To verify the the convergence rate and compare with the nested-loop actor-critic algorithm \cite{yang2019global} which updates the critic more frequently and is usually observed to be less data efficient in practice, we perform a numerical simulation on a small-scale synthetic problem. For the system transition matrices and cost matrices, we choose
\begin{align}
A=\left[\begin{array}{ccc}
0.5 & 0.01 & 0 \\
0.01 & 0.5 & 0.01 \\
0 & 0.01 & 0.5
\end{array}\right], \quad B=\left[\begin{array}{ll}
1 & 0.1 \\
0 & 0.1 \\
0 & 0.1
\end{array}\right], \quad Q=I_{3},\quad R=I_{2}.
\end{align}

We run Algorithm \ref{Alg:AC_LQR} for 2000 iterations, with randomly initialized $K_0$, $\hat{J}_0$ and $\hat{\Omega}_0$ such that the initial $K_0$ is a stabilizing control matrix. The nested-loop actor-critic algorithm that we consider is essentially a nested-loop version of Algorithm \ref{Alg:AC_LQR} where the critic is updated 10 times per actor update with a carefully selected step size to make the estimate of the actor gradient more accurate.

In Figure \ref{fig:lqr_rate}, we plot the convergence rate of the online and nested-loop actor-critic algorithm, where the plot on the left shows the convergence rate of $J(K_k)-J^{\star}$ itself, and on the right, we plot the quantity $\frac{J(K_k)-J^{\star}}{(k+1)^{-2/3}\log^2(k+1)}$. 
The first plot illustrates the superior empirical performance of the online algorithm.
The curve in the second plot approaches a flat line after 500 iterates, which suggests that $J(K_k)-J^{\star}$ indeed roughly decays with rate $\Ocal(\log^2(k+1)/(k+1)^{2/3})$.

\section{Acknowledgement}

This work was partially supported by ARL DCIST CRA W911NF-17-2-0181.


\bibliographystyle{siamplain}
\bibliography{references,ref_Doan}

    

  




\end{document}